\documentclass{amsart}
\usepackage{amsmath}[1996/11/01]
\usepackage{amssymb,amsthm,amsxtra}
\RequirePackage{amssymb}
\RequirePackage{epsfig}
\RequirePackage{ifthen}

\def\[#1\]{\begin{equation}#1\end{equation}}
\makeatletter
\def\beq{%
   \relax\ifmmode
      \@badmath
   \else
      \ifvmode
         \nointerlineskip
         \makebox[.6\linewidth]%
      \fi
      $$
   \fi
}
\def\eeq{%
   \relax\ifmmode
      \ifinner
         \@badmath
      \else
         $$
      \fi
   \else
      \@badmath
   \fi
   \ignorespaces
}

\def\enddisplaymath{\eeq\global\@ignoretrue}
\makeatother

\newtheorem{thm}{Theorem}
\newtheorem{cor}[thm]{Corollary}
\newtheorem{lem}[thm]{Lemma}
\newtheorem{prop}[thm]{Proposition}

\theoremstyle{remark}
\newtheorem*{rem}{Remark}
\newtheorem*{thm-nonum}{Theorem}

\theoremstyle{definition}

\numberwithin{equation}{section}
\numberwithin{thm}{section}
\numberwithin{eg}{section}

\DeclareMathOperator{\Aut}{Aut}

\begin{document}

\title{Quadratic transformations of Macdonald and Koornwinder polynomials}
\author{Eric M. Rains
}
\address{
UC Davis\\
 Department of Mathematics \\
One Shields Ave \\
Davis, CA 95616-8633}
\email{rains@math.ucdavis.edu}
\author{Monica J. Vazirani}
\address{
UC Davis\\
 Department of Mathematics \\
One Shields Ave \\
Davis, CA 95616-8633}
\email{vazirani@math.ucdavis.edu}

\date{June 8, 2006}

\begin{abstract}

When one expands a Schur function in terms of the irreducible characters of
the symplectic (or orthogonal) group, the coefficient of the trivial
character is 0 unless the indexing partition has an appropriate form.  A
number of $q$-analogues of this fact were conjectured in \cite{bcpoly}; the
present paper proves most of those conjectures, as well as some new
identities suggested by the proof technique.  The proof involves showing
that a nonsymmetric version of the relevant integral is annihilated by a
suitable ideal of the affine Hecke algebra, and that any such annihilated
functional satisfies the desired vanishing property.  This does not,
however, give rise to vanishing identities for the standard nonsymmetric
Macdonald and Koornwinder polynomials; we discuss the required modification
to these polynomials to support such results.

\end{abstract}

\maketitle

\tableofcontents



\newcommand{\enbrac}[2]{\langle #1, #2 \rangle}
\newcommand{\<}{\langle}
\renewcommand{\>}{\rangle}
\newcommand{\ssii}[1]{s_{2 #1} s_{2 #1 +1} s_{2 #1 -1} s_{2 #1} }
\newcommand{\ssi}{s_{2 i} s_{2i+1} s_{2i -1} s_{2i} }
\newcommand{\sss}{s_{2} s_{3} s_{1} s_{2} }
\newcommand{\isom}{\simeq}
\newcommand{\WJ}{W^J}
\newcommand{\U}{U}
\newcommand{\Rw}{R}
\newcommand{\RF}{R^F}
\newcommand{\Ri}{R^\iota}
\newcommand{\orbR}{\bar R} 
\newcommand{\wt}{wt}
\newcommand{\hgt}{ht}
\newcommand{\f}{f}
\newcommand{\g}{g}
\newcommand{\dd}{d}
\newcommand{\comment}[1]{\footnote{{\bf ?} #1 {\bf ?}}}
\newcommand{\I}{{\mathcal I}}
\newcommand{\Z}{{\mathbb Z}}
\newcommand{\Hn}{H_n}
\newcommand{\HnC}{H_n^C}
\newcommand{\Hnn}{H_{2n}}
\newcommand{\Hj}[1]{H_{{#1}}}
\newcommand{\Hlam}{H(\lambda)}
\newcommand{\Clam}{\C(q,t)_\lambda}
\newcommand{\C}{{\mathbb C}}
\newcommand{\bT}{\overline T} 
\newcommand{\Ti}{T_i}
\newcommand{\anti}{\tau}
\newcommand{\shift}{\pi}
\newcommand{\Vlam}{V_{<\lambda} }
\newcommand{\Vlesslam}{V_{\le\lambda} }
\newcommand{\Vquot}{\Vlesslam / \Vlam }
\newcommand{\LL}{{\mathcal L}}
\newcommand{\Llam}{{\mathcal L}_{\lambda} }
\newcommand{\vect}{v}
\newcommand{\vw}{\vect_w}
\newcommand{\vwsi}{\vect_{w s_i }}
\newcommand{\vwj}[1]{\vect_{{#1}} }
\newcommand{\eps}[1]{\varepsilon_{#1}}
\newcommand{\epsi}{\eps{i}}
\newcommand{\epsj}{\eps{j}}

\newcommand{\Q}{\mathbb Q}
\newcommand{\R}{\mathbb R}
\newcommand{\F}{\mathbb F}
\newcommand{\T}{\mathbb T}
\renewcommand{\P}{\mathbb P}

\section{Introduction}

Whenever one considers an identity of Schur functions, it is natural to
consider whether that identity admits a $q$-analogue; that is, whether
there is a corresponding identity for Macdonald polynomials.  One such
(classical) identity arises in the representation theory of real Lie
groups, or equivalently in the theory of compact symmetric spaces.

\begin{thm-nonum} \cite{MacD.book}
For any integer $n\ge 0$ and partition $\lambda$ with at most $n$ parts,
the integral
\[
\notag
\int_{O\in O(n)} s_\lambda(O) dO
\]
(with respect to Haar measure on the orthogonal group) vanishes unless
$\lambda=2\mu$ for some $\mu$ (that is, unless every part of $\lambda$ is
even).  Similarly, for $n$ even, the integral
\[
\notag
\int_{S\in Sp(n)} s_{\lambda}(S) dS
\]
vanishes unless $\lambda=\mu^2$ for some $\mu$.
\end{thm-nonum}

Recall the Schur function $s_\lambda$ is a symmetric polynomial in
$n$ variables which gives the character of an irreducible (polynomial)
representation of $U(n)$ 
($GL(n)$)%
.   The character's value on  a matrix is given my evaluating the 
Schur function 
at the matrix's eigenvalues.
The above theorem describes which representations have $O(n)$ ($Sp(n)$)
invariants--exactly those indexed by partitions all of whose parts are
even (occur with even multiplicity).

The symmetric function interpretation of this theorem is that
if one expands $s_\lambda$ in terms of the irreducible characters
of $O(n)$ ($Sp(n)$),
the coefficient of the trivial character is $0$ unless $\lambda = 2 \mu$
($\lambda = \mu^2$).
This formulation has a nice $q$-analogue in several cases.

\begin{rem}
The nonzero values of the integral are in this case all equal to 1; this
will fail upon passing to the Macdonald analogue, although in all cases for
which we can compute the nonzero values, said values are at least ``nice''
(i.e., expressible as a ratio of products of binomials).
\end{rem}

This can in turn be restated in terms of the eigenvalue densities of the
orthogonal and symplectic groups.  For the symplectic group, this is
particularly simple
(integrating over the torus $T$ instead of the whole group):
\[
\int s_\lambda(z_1,z_1^{-1},z_2,z_2^{-1},\dots,z_n,z_n^{-1})
\prod_{1\le i\le n} |z_i-1/z_i|^2
\prod_{1\le i<j\le n} |z_i+1/z_i-z_j-1/z_j|^2
dT
\]
vanishes unless $\lambda=\mu^2$ for some $\mu$.  For the orthogonal group,
the situation is more complicated, as the orthogonal group has two
components, and the structure of the eigenvalues on a given component depends 
significantly on the parity of the dimension; we thus obtain four different
integrals:
\begin{align}
\int s_\lambda(\dots ,z_i^{\pm 1},\dots)
\prod_{1\le i<j\le n}
&
|z_i+1/z_i-z_j-1/z_j|^2
dT
\\
\int s_\lambda(\dots ,z_i^{\pm 1},\dots,\pm 1)
\prod_{1\le i\le n-1}
&
|z_i-1/z_i|^2
\prod_{1\le i<j\le n-1} |z_i+1/z_i-z_j-1/z_j|^2
dT
\\
\int s_\lambda(\dots ,z_i^{\pm 1},\dots,1)
\prod_{1\le i\le n}
&
|1-z_i|^2
\prod_{1\le i<j\le n} |z_i+1/z_i-z_j-1/z_j|^2
dT
\\
\int s_\lambda(\dots ,z_i^{\pm 1},\dots,-1)
\prod_{1\le i\le n/2}
&
|1+z_i|^2
\prod_{1\le i<j\le n} |z_i+1/z_i-z_j-1/z_j|^2
dT
\end{align}
where the first two integrals correspond to the two components of $O(2n)$,
and the last two integral correspond to the two components of $O(2n+1)$,
and the claim is that each integral vanishes unless all ($2n$ or $2n+1$)
parts of $\lambda$ have the same parity.

In \cite{bcpoly}, $q$-analogues of each of these integrals were
conjectured; that is, suitable choices of density were found such that
specializing a Macdonald polynomial as above then integrating against the
appropriate density gives 0 unless the partition satisfies the appropriate
condition.  In particular, the $q$-analogue of the symplectic vanishing
integral (which we will prove in section \ref{sec-sp.vanish}) reads as follows.

\begin{thm-nonum}
\label{thm-Sn}
For any integer $n\ge 0$, and partition $\lambda$ with at most $2n$ parts,
and any complex numbers $q$, $t$ with $|q|,|t|<1$, the integral
\begin{gather*}
\int
P_\lambda(\dots,z_i^{\pm 1},\dots;q,t)
\prod_{1\le i\le n} \frac{(z_i^{\pm 2};q)}
                         {(t z_i^{\pm 2};q)}
\prod_{1\le i<j\le n} \frac{(z_i^{\pm 1} z_j^{\pm 1};q)}
                           {(t z_i^{\pm 1} z_j^{\pm 1};q)}
dT
\end{gather*}
vanishes unless $\lambda=\mu^2$ for some $\mu$.
\end{thm-nonum}

The proof 
 then suggests other statements along these lines, some 
of which are conjectured in  \cite{bcpoly},
 but some of
which are new.

In many of these other identities, we relate a  Macdonald
or Koornwinder polynomial with one value of parameters $q,t$
to polynomials in which $q$ or $t$ is replaced by its
square or square root and thus
these identities can be viewed as ``quadratic" identities
in the sense of basic hypergeometric series.

Similarly, one of the special cases of Theorem \ref{thm-Sn}
proved in \cite{bcpoly} was shown to be equivalent to a 
quadratic transformation for a univariate hypergeometric
series. Thus in a sense these identities can be viewed
as multivariate analogues of quadratic transformations.

There is a fundamental obstruction in using the affine Hecke algebra
approach to directly proving the orthogonal 
cases which here only follow from the observation of \cite{bcpoly} 
that the symplectic and orthogonal identities
 are equivalent by a sort of duality.
Etingof (personal communication) has suggested an
alternate approach using the construction of Macdonald polynomials 
in Etingof-Kirillov \cite{EK}.
This approach works (aside from checking some details)
for the orthogonal but {\it not\/} symplectic vanishing identities,
and
like our approach, it also gives no information about the
nonzero values of the integrals. 
Presumably others of the identities we prove below could be
proved in similar ways, where for Koornwinder polynomials
we must use the 
construction of Oblomkov-Stokman \cite{OS}.

In \cite{bctheta} the Koornwinder polynomials are generalized to
a  family of bi-orthogonal abelian functions.
It is thus natural to conjecture that the vanishing identities
should extend to the elliptic level.
At present, this is somewhat problematic as neither the
(double) affine Hecke algebra approach nor the 
Oblomkov-Stokman construction have been extended to
this setting.

It is also worth noting that a {\it different\/} limit of the bi-orthogonal
abelian functions gives ordinary symmetric Macdonald polynomials
(as orthogonal polynomials) \cite{limits},
suggesting that our Macdonald polynomial identities should
also be limits of elliptic vanishing identities.
It is likely that taking different limits of a single elliptic
vanishing identity could give both a Macdonald and a Koornwinder
identity.
A particularly likely example are the identities
\eqref{Sn1} and \eqref{IK2}.

Identities \eqref{IK} and its dual \eqref{IKdual}
below can be generalized using a third approach that actually works
on both cases.  This will be discussed in a future paper.

Acknowledgments: We would like to thank the Institute for Quantum Information
and the Department of Mathematics at Caltech for hosting our 
respective visits there, where this collaboration began. 
EMR was supported in part by NSF grant DMS-0401387.
MJV was supported in part by NSF grant DMS-0301320, and the UC
Davis Faculty Development Program.

\section{Conventions and Notation}
\label{sec-conventions}


A {\it partition\/} with $\le n$ parts is a nonincreasing integer
tuple $\lambda = (\lambda_1 \ge \lambda_2 \ge \cdots \ge \lambda_n \ge 0)$.
We write $|\lambda| = \sum_{i=1}^n \lambda_i$ 
 or as $\lambda \vdash \sum_{i=1}^n \lambda_i$.
We also let $\ell(\lambda) = \max\{k \ge 0 \mid \lambda_k \neq 0\}$ so
for instance, above we are taking $\ell(\lambda) \le n$.
We will denote the zero (or empty) partition by $0$, when clear in context.
We can picture a partition $\lambda$ as a Ferrer's diagram:
a collection of $|\lambda|$ cells whose coordinates we label
$(i,j)$ with $1 \le j \le \lambda_i$.
So we can refer to a cell as $(i,j) \in \lambda$.
We write $\lambda'$ for the conjugate partition, which
corresponds to a Ferrer's diagram with cells having coordinates
$(j,i)$.

A tuple $\nu = (\nu_1, \ldots, \nu_n)$ of non-negative integers 
is called a {\it composition\/} of $|\nu| = \sum_i \nu_i$.
We will denote by $\nu^+$ the partition obtained by writing the parts
of $\nu$ in nonincreasing order.

Given a partition $\mu$, we write 
$\lambda = \mu^2$
if $\lambda_{2i-1} = \lambda_{2i} = \mu_i$.  
In particular, the parts of $\mu^2$ occur with even multiplicity.
We write $\lambda = 2 \mu$ if $\lambda_i = 2 \mu_i$, so each
part of $2 \mu$ is even.  Note that if $\lambda = \mu^2$ then
the transposed partition $\lambda' = 2 \mu'$.


We define 
$$(a;q) = \prod_{k \ge 0} (1-a q^k)$$
and 
$(a_1, a_2, \ldots, a_\ell ; q) = 
(a_1;q)(a_2;q) \cdots (a_\ell;q).$
As an example, $(x_i^{\pm 1} x_j^{\pm 1}; q) =
(x_i x_j, x_i x_j^{- 1}, x_i^{-1} x_j, x_i^{-1} x_j^{-1}; q) 
=
(x_i x_j; q) (x_i x_j^{- 1}; q) (x_i^{-1} x_j; q) (x_i^{-1} x_j^{-1}; q)$. 
We write $(a;q)$ for what is often denoted $(a;q)_\infty$
in the literature, but as every $q$-symbol we use is infinite,
there is no risk of confusion.

We also define 
\begin{gather*}
C_\mu^0 (x;q,t)= 
\prod_{(i,j) \in \mu} (1-q^{j-1}t^{1-i}x) =
\prod_{1\le i \le \ell(\mu)}\frac{(t^{1-i}x;q)}{(q^{\mu_i}t^{1-i}x;q)},
\\
C_\mu^- (x;q,t)= 
\prod_{(i,j) \in \mu} (1-q^{\mu_i-j}t^{\mu_j'-i}x)
= \prod_{1\le i \le \ell(\mu)}
\frac{(x;q)}{(q^{\mu_i} t^{\ell(\mu) -i} x; q)}
 \prod_{1\le i< j \le \ell(\mu)}
\frac{(q^{\mu_i - \mu_j} t^{j-i}x;q)}{(q^{\mu_i-\mu_j} t^{j-i-1} x; q)},
\\
C_\mu^+ (x;q,t)= 
\prod_{(i,j) \in \mu} (1-q^{\mu_i+j-1}t^{2-\mu_j'-i}x)
= \prod_{1\le i \le \ell(\mu)}
\frac{(q^{\mu_i } t^{2 - \ell(\mu) -i}x;q)}{(q^{2\mu_i} t^{2-2i} x; q)}
 \prod_{1\le i< j \le \ell(\mu)}
\frac{(q^{\mu_i +\mu_j} t^{3-j-i}x;q)}{(q^{\mu_i+\mu_j} t^{2-j-i} x; q)}.
\end{gather*}
Similar to the $q$-symbols, we let
$C_\mu^{0,\pm}
(a_1, a_2, \ldots, a_\ell ; q) = 
C_\mu^{0,\pm}(a_1;q)C_\mu^{0,\pm}(a_2;q) \cdots C_\mu^{0,\pm}(a_\ell;q).$
We refer to \cite{bcpoly} for more details about these expressions
and relations that hold among them (in particular those expressing
$C_{\mu^2}^{0,\pm}(x;q,t) $ or $C_{2 \mu}^{0,\pm}(x;q,t)$ 
in terms of $C_{\mu}^{0,\pm}(x;q^2,t) $ and  $C_{\mu}^{0,\pm}(x;q,t^2)$).

It will be convenient in the sequel to use a plethystic substitution
notation slightly different from that in the literature.
When  we write $g([r_k])$ for symmetric functions $g, r_k, k \ge 1$
we mean the image of $g$ under the homomorphism $p_k \mapsto r_k$
where the $p_k$ are the power sum symmetric functions. 
We take the convention $p_x = 0$ if $x \not\in \{ 1, 2, 3, \ldots \}$.
We abbreviate the case $r_{2k+1} = 0$, $r_{2k} = p_k$ by
$g([2 p_{k/2}])$.
This is plethystic notation for the specialization
$g(\cdots , \pm \sqrt{x_i}, \cdots)$.

\subsection{The extended affine symmetric group $\widetilde{S}_n$}
$\widetilde{S}_n=S_n\ltimes\Z^n$ can be identified with the group of
bijections $w:\Z\to\Z$ such that $w(x+n)=w(x)+n$ for all
$x\in \Z$; if we also include bijections such that $w(x+n)=w(x)-n$,
we obtain a group $\widetilde{S}_n^+=(S_n\times \Z_2)\ltimes \Z^n$,
which is also an extended affine Weyl group (see section \ref{sec-extended}).
The length 0 subgroup of  $\widetilde{S}_n^+$
 is generated by $\shift(x)=x+1$ and $\iota(x)=n+1-x$.

$\widetilde{S}_n$ has generators $s_0, s_1, \ldots, s_{n-1}, \shift$ where
\begin{gather*}
s_j(i) = \begin{cases} i & \quad  i \not\equiv j, j+1 \mod n \\
		i+1 & \quad   i \equiv j \mod n\\
		i-1 &  \quad i \equiv j+1 \mod n.\\
\end{cases}
\end{gather*}
By convention we will view these bijections as acting on $\Z$
from the {\it right}.

It is easy to see these generators satisfy the type $A$ braid relations
\begin{gather*}
s_i s_j = s_j s_i \qquad i-j \not\equiv \pm 1\pmod  n \\
s_i s_j s_i = s_j s_i s_j \qquad i-j \equiv \pm 1\pmod  n; n > 2,
\end{gather*}
and
\begin{gather*}
\shift s_i \shift^{-1} = s_{i-1}.
\end{gather*}
and quadratic relation $s_i^2 = 1$.

\label{sec-define.hecke}
The {\it extended affine Hecke algebra $\Hn$  of type $A$} is defined
to be the $\C(q,t)$-algebra with
generators $T_0$,$T_1$,\dots $T_{n-1}$,$\shift$, subject to the braid
relations
\begin{gather*}
T_i T_j = T_j T_i \qquad i-j \not\equiv \pm 1\pmod  n \\
T_i T_j T_i = T_j T_i T_j \qquad i-j \equiv \pm 1\pmod  n; n > 2,
\end{gather*}
and the
quadratic relation
\begin{gather*}
(T_i - t) (T_i + 1) = 0,
\end{gather*}
 and
\[
\shift T_i \shift^{-1} = T_{i-1}.
\]

Given a reduced word $u = s_{i_1} s_{i_2} \cdots s_{i_k}$,
we write $T_u =  T_{i_1} T_{i_2} \cdots T_{i_k}$,
which is independent of reduced word expression by the relations above.
Note $T_u T_v = T_{uv}$ if $\ell( u) + \ell(v)=\ell(uv)$.

Observe that on specializing $t=1$ we recover the group algebra
$\C(q) \widetilde S_n$ whose generators we typically denote
$\{s_0, s_1, \ldots, s_{n-1}, \shift\}$.

Given an automorphism $\phi: \Hn \to \Hn$ and right module $\LL$,
we write $\LL^\phi$ for the twisted module with action
$v \cdot h := v( \phi(h) )$.
In the case $\phi(h) = T_u h T_u^{-1}$
we write $\LL^u$ for $\LL^\phi$.

We will write $\bT_i = T_i + 1 -t$. Note $\bT_i T_i = t$.
(This is {\it not\/} Lusztig's bar involution.)

We have another presentation of  $\Hn$ given by
$T_1, T_2, \ldots, T_{n-1}$, $Y_1^{\pm 1}, Y_2^{\pm 1}, \ldots,  Y_n^{\pm 1}$ 
with additional relations
\begin{align*}
Y_i Y_j = Y_j Y_i  &\qquad \forall i, j \\
T_i Y_j = Y_j T_i &\qquad j \neq i, i+1 \\
T_i Y_i^{-1} T_i = t Y_{i+1}^{-1} &\qquad 1 \le i < n, 
\end{align*}
where we can also express the final one as $T_i Y_{i+1} = Y_i \bT_i$.

This presentation relates to the first via:
\begin{align}
\label{eq-Y1}
Y_1&=T_1\dots T_{n-1}\shift,\\
\label{eq-Y2}
Y_2&=T_2\dots T_{n-1}\shift \bT_1,\\
\notag
&\vdots\\
\label{eq-Yn}
Y_n&=\shift\bT_1\dots \bT_{n-1}.
\end{align}
(This disagrees with the convention that $tY_{i+1}=T_i Y_i T_i$, but has
the advantage of making dominant weights map to positive words!) 
That is, for a partition $\lambda$, $Y^\lambda = Y_1^{\lambda_1}
Y_2^{\lambda_2} \cdots  Y_n^{\lambda_n}$ simplifies in the other
generators to a word involving only $T_i$ and $\shift$  and not
involving $\bT_i$.

$\Hn $ acts on the space of polynomials
$V= \C(q^{1/n},t)[x_1, \ldots, x_n, (x_1 x_2 \cdots x_n)^{-1/n}]$
via:
\begin{align}
\label{eq-Ti}
T_i f &= tf + \frac{x_{i+1}-t x_i}{x_{i+1}-x_i} (f^{s_i}-f)\\
T_0 f &= tf + \frac{x_{1}-t q x_n}{x_{1}- q x_n} (f^{s_0}-f) \\
(\shift f)(x_1,\dots x_n) &= f(q x_n,x_1,\dots x_{n-1}),
\end{align}
where $f^{s_i}(x_1, \ldots, x_i, x_{i+1}, \ldots, x_n) =
f(x_1, \ldots,  x_{i+1}, x_i,\ldots, x_n)$
and $f^{s_0}(x_1, \ldots,  x_n) =
f(q x_n, \ldots, q^{-1} x_1)$.
Observe
\begin{align}
T_i 1 &= t  \\
\shift 1 &= 1 \\
Y_i 1 &= t^{n-i} .
\end{align}
Observe the $T_i$ act trivially on $(x_1 x_2 \cdots x_n)^{-1/n}$,
but $\shift$ multiplies it by $q^{-1/n}$.

Given a partition $\lambda$ or more generally a dominant weight of
$SL_n \times GL_1$,
i.e. a nonincreasing sequence of rational numbers with integer sum
and integer differences, we can associate a monomial in $V$, namely
$\prod_i x_i^{\lambda_i}$.
This generates a $\widetilde S_n$-submodule of $V$.
This however is not invariant under $\Hn$, but if we sum the spaces
corresponding to all weights weakly dominated by $\lambda$
then the space is invariant under $\Hn$ and affords a filtration.
The associated graded of the filtration gives a deformation
to $\Hn$ of the $\widetilde S_n$-submodule associated to $\lambda$.
In this space, the commuting operators $Y_i$ have joint
eigenvalues which are simply permutations of
the sequence
\[
\dots q^{\lambda_i} t^{n-i}\dots
\]
Generically, this deformation is a submodule of $V$, and thus the
corresponding eigenfunctions are polynomials, namely the nonsymmetric
Macdonald polynomials.

Recall if $L$ is a functional on the polynomial space (or on any left module),
then $h \in \Hn$ acts on the
right via $(L \cdot h)(f) := L(hf)$.

This representation has the following interpretation in terms of the
double affine Hecke algebra (while we do not use or even define
the double affine Hecke algebra here, it is worth
noting we can view these problems in a larger context).
 Our affine Hecke algebra is a 
subalgebra of the double affine Hecke algebra, 
and it
has 
a ``trivial" module, which
 is the one dimensional module on which
$\shift- 1$ and all $T_i - t$ vanish.
If we induce this trivial module up to the
double affine Hecke algebra and then restrict it back down,
$V$ sits inside the restriction.
The Mackey formula thus gives us  a decomposition of  $V$
into irreducibles (when $q$, $t$ are generic) which we describe
explicitly below.


If we specialize $q$, $t$ to complex numbers such that $|q|, |t| < 1$,
then the nonsymmetric density
\[
\Delta_S = 
\Delta_S^{(n)}(q,t) = 
\prod_{i<j}
\frac{(x_i/x_j,qx_j/x_i;q)}{(tx_i/x_j,qtx_j/x_i;q)}
\]
is defined and can be integrated over the unit torus.
Moreover it is a standard result of Macdonald  polynomials
theory that if $f \in \C(q^{1/n},t)[x_1, \ldots, x_n,  (x_1 \cdots x_n)^{-1/n}]$
then
\begin{gather}
\frac{ \int f \Delta_S dT}{\int \Delta_S dT} \in \C(q^{1/n},t).
\end{gather}
That is, there exists a rational function in $q^{1/n}, t$
that agrees with the above for any specialization such that the
integrals are defined. 
Similar comments apply to all the integrals we consider which
can thus be considered either as analytic quantities with
appropriately specialized parameters or as algebraic quantities
with generic parameters.
In particular, the normalized integral
$\frac{ \int f \Delta_S dT}{\int \Delta_S dT} = [E_0]f$
where $E_0$ is the nonsymmetric Macdonald polynomial corresponding
to the empty partition.
A similar statement holds in the other cases.

Above, we used the notation
$$[f_\mu] g$$
for the coefficient of $f_\mu$ in the expansion of $g$,
where $\{ f_\mu \}$
is a given a basis
of some space of functions and $g$ is  another 
function in that space.
It should be clear in all cases in which we use this notation which
basis is intended.

Note
that $\shift$ is self-adjoint
and the $T_i$ are adjoint to $T_{n-i}$
with respect to
the inner product this density defines:
$$\langle f, g \rangle =
\int f(x_1, \ldots, x_n) g(\frac {1}{x_n}, \ldots, \frac{1}{x_1})
\Delta_S dT,
$$
where $dT$ is Haar measure on the unit torus.

An equivalent way of stating this uses the fact that $\Hn \otimes \Hn$
has a natural action on 
$\C(q^{1/n},t)[y_1, \ldots, y_n, z_1, \ldots,$
$z_n, (y_1 z_1 \cdots y_n z_n)^{-1/n}]$
and says that
the linear functional
\begin{gather*}
L(h) = \int h(x_1, \ldots, x_n, \frac {1}{x_n}, \ldots, \frac{1}{x_1})
\Delta_S^{(n)}(q,t) dT
\end{gather*}
is annihilated by the ideal $\langle
\shift \otimes 1 - 1 \otimes \shift, \,
T_i \otimes 1 - 1 \otimes T_{n-i}, (1 \le i\le n)
\rangle$.
As we will see below, such annihilation gives rise to 
vanishing identities. 
  In this case, we obtain the (standard) fact  that for
weights 
$\lambda$ and $\mu$
\begin{gather*}
\int P_\lambda(x_1, \ldots, x_n) P_\mu( \frac {1}{x_n}, \ldots, \frac{1}{x_1})
\widetilde \Delta_S dT
\end{gather*}
vanishes if $\lambda \neq \mu$.
Here, $\widetilde \Delta_S$ is the symmetric density
\[
\widetilde \Delta_S=
\widetilde \Delta_S^{(n)}(q,t)=
\prod_{i<j}
\frac{(x_i/x_j,x_j/x_i;q)}{(tx_i/x_j,tx_j/x_i;q)}
= \prod_{i \neq j} \frac{(x_i/x_j;q)}{(tx_i/x_j;q)}
\]
which up to scalar is the symmetrization of $\Delta_S$.

The operators $Y_i$ are {\it not\/} self-adjoint and more generally the
ideal does not contain elements of the form 
$Y_i \otimes 1 - c 1 \otimes Y_j^{\pm 1}$.
However, if we conjugate by $1 \otimes T_{w_0}$ 
 it will contain
$Y_i \otimes 1 -  1 \otimes Y_i$ and this implies orthogonality
of $Y$-eigenvectors 
with respect to the conjugated inner product.
With respect to the original inner product, we find that the
eigenfunctions of the $Y_i$ are orthogonal to
the images of those functions under $T_{w_0}^{- 1}$.
This is precisely the orthogonality of nonsymmetric Macdonald
polynomials given in \cite{Chered}.
(To be precise, Cherednik shows that the nonsymmetric Macdonald polynomials 
(a.k.a. the eigenfunctions of the $Y_i$) are
orthogonal to the polynomials modified 
by the substitution $q \to \frac 1q$, $t \to \frac 1t$,
but this turns out to be equivalent.)
It follows that the symmetric Macdonald polynomials are
orthogonal with respect to this density and hence the
symmetrized density.

\subsection{The extended affine hyperoctahedral group $\widetilde{C}_n$}

We also consider $\widetilde{C}_n=C_n\ltimes \Z^n$, which can be
identified with the centralizer in $\widetilde{S}_{2n}$ of the element $\iota$
of $\widetilde{S}_{2n}^+$, or equivalently
as the group of
bijections $w: \Z \to \Z$ such that $w(i+2n) = w(i) + 2n$
and $w(2n+1-i) = 2n+1-w(i)$.
$\widetilde{C}_n$ has generators $s_0, s_1, \ldots, s_{n},$ where
\begin{gather*}
s_j(i) = \begin{cases} i & \quad i \not\equiv j,j+1,2n-j,2n+1-j \mod 2n \\
		i+1 & \quad  i \equiv j,2n-j \mod 2n\\
		i-1 & \quad  i \equiv j+1,2n+1-j \mod 2n.\\
\end{cases}
\end{gather*}
It is easy to see these generators satisfy the type $C$ braid relations
\begin{align*}
& s_i s_j = s_j s_i  \qquad \qquad  |i-j| > 1  \\
& s_i s_j s_i = s_j s_i s_j  \qquad |i-j| = 1, i, j \neq 0,n  \\
& s_0 s_1 s_0 s_1 = s_1 s_0 s_1 s_0     \\
& s_n s_{n-1} s_n s_{n-1} = s_{n-1} s_n s_{n-1} s_n   
\end{align*}
and quadratic relation $s_i^2 = 1$.

For $n>1$, the
{\it
affine Hecke algebra $\HnC$  of type $BC$}
is defined to be the $\C(q,t, a, b,c,d)$-algebra with
generators $T_0$,$T_1$,\dots $T_{n}$, 
subject to the type C
braid relations
\begin{gather*}
T_i T_j = T_j T_i \qquad |i-j| > 1 \\
T_i T_j T_i = T_j T_i T_j \qquad |i-j| = 1, i, j \neq 0,n \\
T_0 T_1 T_0 T_1 = T_1 T_0 T_1 T_0 \\
T_n T_{n-1} T_n T_{n-1} = T_{n-1} T_n T_{n-1} T_n 
\end{gather*}
In fact, the algebra $\HnC$ can be defined for $n=1$ and all considerations
below will work in that case.  If $n=1$ there are simply no braid
relations, only quadratic ones.  We omit the details.

The diagram automorphism of affine $C_n$ gives rise to an
action of the involution $\sigma$ on $\HnC$ given by
\begin{gather*}
\sigma T_i \sigma^{-1} = T_{n-i} \\
\sigma a \sigma^{-1} =  c/\sqrt{q} \\
\sigma b \sigma^{-1} =  d/\sqrt{q}
\end{gather*}
If the action of $\sigma$ on scalars is trivial,
i.e. $c=a\sqrt{q}$, $d=b\sqrt{q}$,
then we can enlarge $\HnC$ to an {\it extended\/} affine Hecke
algebra  as in section \ref{sec-extended}.
In general, we can view $\sigma$ as giving an intertwiner
between Hecke algebras with different parameters. 

For $1\le i < n$, we will write $\bT_i = T_i + 1 -t$. Note $\bT_i T_i = t$.
We set $\bT_n = T_n + 1 + ab$ so that $\bT_n T_n = -ab$ and
$\bT_0 = T_0 + 1 + \frac {cd}{q}$ so $\bT_0 T_0 = - \frac {cd}{q}$.

As with type $A$, 
we have another presentation of  $\HnC$ given by
$T_1, T_2, \ldots, T_{n}$, $Y_1^{\pm 1}, Y_2^{\pm 1}, \ldots,  Y_n^{\pm 1}$ 
with additional relations
\begin{gather*}
Y_i Y_j = Y_j Y_i \qquad \forall i, j \\
T_i Y_j = Y_j T_i \qquad j \neq i, i+1 \\
T_i Y_i^{-1} T_i = t Y_{i+1}^{-1} \qquad 1\le i < n \\
T_n Y_n^{-1} T_n = - t^{2-2n} \frac {q}{cd} Y_n-
t^{1-n}(\frac {q}{cd} +1) T_n.
\end{gather*}

This presentation relates to the first via:
\begin{align}
Y_1 &= T_1 \dots T_n\dots T_1 T_0\\
Y_2 &= T_2 \dots T_n\dots T_1 T_0\bT_1\\
\notag
& \; \vdots\\
Y_n &= T_n \dots T_1 T_0 \bT_1\dots \bT_{n-1}.
\end{align}

  There is also an intertwiner
\[
Y_{\omega} = 
\prod_{1\le i\le n} T_n \dots T_i
\sigma
\]
which satisfies
$Y_i Y_\omega = Y_\omega Y_i$ (but note that the two $Y_i$
live in different Hecke algebras) and
$Y_{\omega}^2=Y_1 Y_2 \dots Y_n$.
The significance of the intertwiner $Y_\omega$ is that the
symmetric version: $Y_\omega (1+Y_1^{-1})(1+Y_2^{-1}) \cdots (1+Y_n^{-1})$
takes a Koornwinder polynomial with parameters $a,b,c,d$ to
a multiple of the corresponding Koornwinder polynomial with
parameters $c/\sqrt{q}, d/\sqrt{q},a\sqrt{q},b\sqrt{q}$.
In fact, this is precisely the difference operator which was
the fundamental tool of \cite{bcpoly}. 
This is one reason why in section \ref{sec-extended} we
consider such a general version of extended affine Weyl groups.

When computing in  $\HnC$ or more generally in the braid group,
$B(\widetilde{C}_n)$, one helpful tool is
the natural injection $\HnC \to \Hnn$ such that
\begin{gather*}
\quad T_i\mapsto T_i T_{2n-i} \\
T_0\mapsto T_0 \\
T_n\mapsto T_n
\end{gather*}
and such that $\sigma$ acts as conjugation by $\shift^n$.  Under this
mapping, the natural lifting of the $Y$ operators to the braid group
behaves as follows:
\begin{align}
Y_i&\mapsto Y_i Y_{2n-i}^{-1}\\
Y_\omega\sigma^{-1}&\mapsto \prod_{1\le i\le n} Y_i \shift^{-n}.
\end{align}

The  Hecke algebra $\HnC$ and the intertwiner $\sigma$ 
act on Laurent  polynomials 
$\C(q^{1/2},t,a,b,c,d)[x_1^{\pm 1}, \ldots, x_n^{\pm n}] $ via:
\begin{align}
T_0 f &= -(cd/q)f+\frac{(1-c/x_1)(1-d/x_1)}{1-q/x_1^2} (f^{s_0}-f)\\
T_i f &= tf + \frac{x_{i+1}-t x_i}{x_{i+1}-x_i} (f^{s_i}-f)\\
T_n f &= -abf+\frac{(1-a x_n)(1-b x_n)}{1-x_n^2} (f^{s_n}-f)\\
 (\sigma f)(a,b,c,d; x_1,\dots x_n) &=
f(c/\sqrt{q}, d/\sqrt{q},a\sqrt{q},b\sqrt{q};\sqrt{q}/x_n,\dots \sqrt{q}/x_1) 
\end{align}
Recall
$f^{s_i}(x_1, \ldots, x_i, x_{i+1}, \ldots, x_n) =
f(x_1, \ldots,  x_{i+1}, x_i,\ldots, x_n)$ for $i \neq 0,n$,
$f^{s_0}(x_1, \ldots,  x_n) = f(q/ x_1, x_2 \ldots,  x_n)$
and $f^{s_n}(x_1, \ldots,  x_n) =  f( x_1, \ldots,x_{n-1}, 1/ x_n)$.

In particular, in the space corresponding to monomials for the partition
$\lambda$, the joint eigenvalues of the operators
$Y_i,(abcdt^{2n-2}q^{-1})Y_i^{-1}$ are (signed) permutations of the sequence
\[
\dots q^{\lambda_i} t^{2n-1-i} (abcd/q)\dots q^{-\lambda_i} t^{i-1}
\]

The nonsymmetric density is
\[
\Delta_K = \Delta_K^{(n)}(a,b,c,d;q,t) =
\prod_{1\le i\le n} \frac{(x_i^2,qx_i^{-2};q)}{(a x_i,b x_i,c x_i,d
  x_i,aqx_i^{-1},bqx_i^{-1},cx_i^{-1},dx_i^{-1};q)} \prod_{1\le i<j\le n}
\frac{(x_ix_j^{\pm 1},q x_i^{-1} x_j^{\pm 1};q)}{(tx_ix_j^{\pm 1},
qt x_i^{-1} x_j^{\pm 1};q)}
\]
Here we integrate over the unit torus with parameters  specialized
to have norm $< 1$.
As with the $S_n$ case, the normalized integral of any polynomial
with rational function coefficients meromorphically continues
to a rational function.

The operators $T_i$ are self-adjoint with respect to the
induced inner product.

The corresponding symmetric density is
\[
\widetilde \Delta_K
=
\widetilde \Delta_K^{(n)}(a,b,c,d;q,t)
=
\prod_{1\le i\le n} \frac{(x_i^{\pm 2};q)}{(a x_i^{\pm 1} ,b x_i^{\pm 1} ,
c x_i^{\pm 1} , d x_i^{\pm 1} ;q)} \prod_{1\le i<j\le n}
\frac{(x_i^{\pm 1} x_j^{\pm 1};q)}{(tx_i^{\pm 1} x_j^{\pm 1}; q)}
\]

Again normalized integrals over this density can be taken by
computing the constant coefficient in the expansion with respect
to Koornwinder polynomials, giving a rational function in $q,t,a,b,c,d$.
In fact, one of the mail results of \cite{bcpoly}
is that this integral essentially depends algebraically on $n$.
More precisely, it is shown there that there exists a functional
$I_K(f; q,t, T; a,b,c,d)$ on the space of ordinary symmetric functions
such that
$$I_K(f; q,t, t^n; a,b,c,d) = [K_0^{(n)}(\cdots z_i \cdots; q,t; a,b,c,d)]
f(\cdots z_i^{\pm 1} \cdots) = \frac {\int f \widetilde \Delta_K}{\int
  \widetilde \Delta_K}$$ for all integers $n \ge 0$.  This can also be
viewed as taking coefficients with respect to a basis $\widetilde
K_\lambda( \, ; q,t,T; a,b,c,d)$ of the space of symmetric functions over
$\C(q,t,T,a,b,c,d)$ with the property that for all integers $n$ such that
$n\ge \ell(\lambda)$, $\widetilde K_\lambda(\cdots z_i^{\pm 1} \cdots;
q,t,t^n; a,b,c,d) =K_\lambda^{(n)}(\cdots z_i \cdots; q,t; a,b,c,d).$
(These $\widetilde K_\lambda$ transform nicely under an analogue of
Macdonald's involution and so we can use them to prove dual results in
several cases.)

\section{A $U(2n)/Sp(2n)$ vanishing integral}
\label{sec-sp.vanish}

\begin{thm}
\label{thm-vanish1}
For any integer $n\ge 0$, and partition $\lambda$ with at most $2n$ parts,
and any complex numbers $q$, $t$ with $|q|,|t|<1$, the integral
\begin{multline}
\int
P_\lambda^{(2n)}(z_1^{\pm 1},\dots, z_n^{\pm 1};q,t) \widetilde \Delta_K^{(n)}
(\sqrt{t}, -\sqrt{t}, \sqrt{qt}, -\sqrt{qt}; q, t) dT
\\
=
\int
P_\lambda^{(2n)}(z_1^{\pm 1},\dots, z_n^{\pm 1};q,t)
\prod_{1\le i\le n} \frac{(z_i^{\pm 2};q)}
                         {(t z_i^{\pm 2};q)}
\prod_{1\le i<j\le n} \frac{(z_i^{\pm 1} z_j^{\pm 1};q)}
                           {(t z_i^{\pm 1} z_j^{\pm 1};q)}
dT
\end{multline}
vanishes unless $\lambda=\mu^2$ for some $\mu$.
\end{thm}

\begin{proof}
Consider the following linear functional on the space of polynomials in
$2n$ variables:
\begin{multline*}
L(f)
:=
\int f(z_1,z_2,\dots,z_n,1/z_n,\dots,1/z_2,1/z_1)
\prod_{1\le i\le n} \frac{(z_i^2,q z_i^{-2};q)}{(t z_i^2,q t z_i^{-2};q)}
\prod_{1\le i<j\le n} \frac{(z_i z_j^{\pm 1},q z_i^{-1} z_j^{\pm 1};q)}
                           {(t z_i z_j^{\pm 1},qt z_i^{-1} z_j^{\pm 1};q)}
dT
\\
=
\int f(z_1,z_2,\dots,z_n,1/z_n,\dots,1/z_2,1/z_1)\Delta_K^{(n)}
(\sqrt{t}, -\sqrt{t}, \sqrt{qt}, -\sqrt{qt}; q, t)
dT
\end{multline*}
If $f$ is symmetric, then we can freely symmetrize the density; since the
density is recognizable as a special case of the nonsymmetric Koornwinder
density, it symmetrizes to the symmetric density above.  In other words, it
will suffice to show that $L(P_\lambda(\, ;q,t))=0$ unless $\lambda=\mu^2$.


Since nonsymmetric Macdonald polynomials of type $C$ are
orthogonal with respect to the density
$\Delta_K^{(n)} \left(\sqrt{t}, -\sqrt{t}, \right.$
$\left. \sqrt{qt}, -\sqrt{qt}; q, t \right)$,
we can interpret the result as saying when we expand
type $A$ Macdonald polynomials in terms of those of type $C$
the coefficient of the trivial one is zero unless $\lambda=\mu^2$.
(In the notation of section \ref{sec-conventions},
$[E_0^C ] E_\lambda^A = 0$ unless $\lambda=\mu^2$.)

The advantage of passing to this nonsymmetric functional is that we can
use the affine Hecke algebra.  Indeed, a straightforward calculation gives
the following facts about the interaction between $L$ and the Hecke
algebra:
\begin{align}
\label{eq-L1}
L(T_0 f) &= t L(f)\\
\label{eq-L2}
L(T_n f) &= t L(f)\\
\label{eq-L3}
L(T_i f) &= L(T_{2n-i} f),\quad 1\le i\le n-1.
\end{align}
But in fact, for generic $q$ and $t$, any linear functional satisfying
these three conditions will also satisfy the vanishing property
``$L(P_\lambda(;q,t))=0$ unless $\lambda=\mu^2$''.

The calculation to verify \eqref{eq-L1}, \eqref{eq-L2}, \eqref{eq-L3}
is very similar to
that of computing the adjoint of $T_i$ with respect to
$\langle \, , \, \rangle$.
A sample computation is given here:
First recall  $(T_n - t)f = \frac{x_{n+1}-t x_n}{x_{n+1}-x_n} (f^{s_n}-f)$
by \eqref{eq-Ti}.
After specializing as above, this will become
$ \frac{z_{n}^{-1} -t z_n}{z_{n}^{-1}-z_n} g_1$
where $g_1$ is  a Laurent polynomial sent to $-g_1$ 
under the change of variables $z_n \leftrightarrow z_n^{-1}$.
Observe the density
\begin{gather*}
\Delta = 
\prod_{1\le i\le n} \frac{(z_i^2,q z_i^{-2};q)}{(t z_i^2,q t z_i^{-2};q)}
\prod_{1\le i<j\le n} \frac{(z_i z_j^{\pm 1},q z_i^{-1} z_j^{\pm 1};q)}
                           {(t z_i z_j^{\pm 1},qt z_i^{-1} z_j^{\pm 1};q)}
= \frac{(1-z_n^2)}{(1-tz_n^2)} g_2
\end{gather*}
where $g_2$ is symmetric under
the change of variables $z_n \leftrightarrow z_n^{-1}$.
Hence $L( (T_n - t)f ) = \int \frac{z_{n}^{-1} -t z_n}{z_{n}^{-1}-z_n} g_1
 \frac{(1-z_n^2)}{(1-tz_n^2)} g_2
= \int g_1 g_2 = \int (-g_1) g_2$
as we are integrating over the torus $T$ and hence get the same
integral 
under the change of variables $z_n \leftrightarrow z_n^{-1}$.
This shows $L( (T_n - t)f ) = 0$.

Let $V_{\le \lambda}$ be the space of polynomials spanned by monomials
$x^\nu:=x_1^{\nu_1}x_2^{\nu_2}\dots$ where $\nu$ is a composition 
of
$|\lambda|$ dominated by $\lambda$ (i.e., such that the corresponding
partition $\nu^+$ is dominated by $\lambda$); similarly let $V_{<\lambda}$
be the
space spanned by monomials {\em strictly} dominated by $\lambda$.  Both
subspaces are invariant under the action of the Hecke algebra, and we may
thus consider the spaces $\LL_\lambda$ 
of functionals on
$V_{\le\lambda}/V_{<\lambda}$
satisfying \eqref{eq-L1},\eqref{eq-L2},\eqref{eq-L3}.
If we can show that  $\LL_\lambda=0$ 
unless $\lambda=\mu^2$, we will be done,
since $V_{\le\lambda}/V_{<\lambda}$ is isomorphic (for generic $q$, $t$) to
the invariant subspace generated by $P_{\lambda}(;q,t)$
(with basis given by  $E_\nu(;q,t)$ with $\nu^+ = \lambda$).

Fix a partition $\lambda$ {\it not\/} of the form $\mu^2$.
Now, the monomials in the orbit of $x^\lambda$ form a basis of
$V_{\le\lambda}/V_{<\lambda}$ (for all nonzero $q, t$, not just generic
$q,t$), and in that basis, the action of the Hecke
algebra has coefficients in
$\Z[q^{\pm 1},t]$. 
Thus $\LL_\lambda$  
is the
solution space of a system of linear equations with coefficients in
$\Z[q^{\pm 1},t]$. 
 Now, the generic dimension of such a solution space is
bounded above by the dimension under any specialization. Therefore, it will
suffice to find   {\it one\/} 
 such specialization for which the claim holds, and thus
the dimension of the solution space is $0$ for generic $q,t$.

In particular, take $t=1$, so that the affine Hecke algebra is just the
group algebra of $\widetilde S_n$
 with the corresponding action on polynomials.  Then any
functional $L\in  \LL_\lambda$ 
is invariant under the subgroup generated by
$s_0$, $s_n$, and $s_i s_{2n-i}$ for $1\le i\le n-1$.  This is precisely
the subgroup of elements invariant under the involution $s_i\mapsto
s_{2n-i}$, and thus in particular contains a number of translations, which
act diagonally on monomials.  It follows immediately that $L(x^\nu)=0$
unless $\nu_i=\nu_{2n+1-i}$ for $1\le i\le n$.  Since $\nu$ is simply a
permutation of $\lambda$, the claim follows.
\end{proof}

\begin{rem}
It similarly follows that for generic $q$, $t$, $\dim(\LL_{\mu^2})\le 1$.  In
fact, since the integrals
\[
\int m_{\mu^2}(\dots z_i^{\pm 1}\dots) dT
\]
are nonzero, we can also conclude that $\dim(\LL_{\mu^2})\ge 1$
for {\em all} $q,t$
(since we have exhibited a linear functional that specializes to a nonzero
functional.)
This implies for a wide class of irreducible representations
we have a multiplicity one condition, i.e., that
there exists at most a $1-$dimensional space of linear
functionals satisfying the above invariance
conditions \eqref{eq-L1},\eqref{eq-L2},\eqref{eq-L3}.
This is in a sense a deformation of the fact that $(U(2n), Sp(2n) )$
is a Gelfand pair,
together with the identification of which representations are
spherical.
This appears to be true for general representations, but
we do not consider that question here.
\end{rem}

Now, as it stands, this argument is somewhat unsatisfactory; it would be
nice to avoid the step of specialization to $t=1$.  Certainly, there is a
natural analogue of the subgroup of translations inside the affine Hecke
algebra; unfortunately, the conditions on $L$ do {\em not} imply any sort of
invariance with respect to the standard commutative subalgebra.  Related to
this is the fact that the obvious corresponding statement for nonsymmetric
Macdonald polynomials does not hold; that is, for $t=1$, the conditions on
$L$ suffice to make $L(E_\nu(;q,1))=0$ unless $\nu_i=\nu_{2n+1-i}$ for
$1\le i\le n$, but this is {\em not} true for $t\ne 1$.  The key to
resolving both of these issues is the fact that, although the standard
commutative subalgebra is in some sense canonical (or, at least, is one of
two canonical choices), it is not the only reasonable choice; we will
consider this in more detail in the sections below.
Equivalently, we can leave the commutative subalgebra alone and
transform the functional, thus conjugating the ideal of equations
on the functional.
This gives nice identities for nonsymmetric Macdonald polynomials
but makes the resulting functional extremely complicated.
We consider this approach in sections \ref{sec-hecke} and \ref{sec-comm}.

Another flaw, which is intrinsic in the way we use the affine Hecke
algebra, is that we obtain no information about the nonzero values of the
integral.  Indeed, the conditions on $L$ determine it only up to a scalar
multiple for each $\mu$.  In this case, the nonzero values were already
determined (conditional on the vanishing result) in \cite{bcpoly}; but for
some of the other vanishing results we prove below, it is still an open
question to determine the nonzero values.

For the present case, however, we have (in the notation of \cite{bcpoly},
and recalling the argument given there)

\begin{cor}
\label{cor-nonzeroSp}
For any integer $n\ge 0$ and any partition $\mu$ with at most $n$ parts,
\begin{multline}
\frac{1}{Z}
\int
P_{\mu^2}(\dots,z_i^{\pm 1},\dots)
\prod_{1\le i\le n} \frac{(z_i^{\pm 2};q)}{(t z_i^{\pm 2};q)}
\prod_{1\le i<j\le n} \frac{(z_i^{\pm 1} z_j^{\pm 1};q)}
                           {(t z_i^{\pm 1} z_j^{\pm 1};q)}
dT
\\
=
\frac{1}{Z}
\int
P_{\mu^2}(\dots,z_i^{\pm 1},\dots)
 \widetilde \Delta_K^{(n)}(\sqrt t, -\sqrt t, \sqrt{qt}, -\sqrt{qt};
q,t)
dT
=
\frac{C^0_\mu(t^{2n};q,t^2)C^-_\mu(qt;q,t^2)}
     {C^0_\mu(qt^{2n-1};q,t^2)C^-_\mu(t^2;q,t^2)},
\end{multline}
where the normalization $Z$ is chosen to make the integral $1$ when
$\mu=0$.
\end{cor}

\begin{proof}
Let $L$ be the above linear functional on symmetric functions (so that we
are computing the value of $L(P_{\mu^2}(\, ;q,t))$); we have already
established that $L(P_\lambda(\, ;q,t))=0$ unless $\lambda$ is of the form
$\mu^2$.  Now, consider the value
\[
L((e_1-e_{2n-1})P_\lambda(\, ;q,t)),
\]
where $e_i$ is the elementary symmetric function.
On the one hand, this is 0, since the mere act of specializing to $z_i^{\pm
  1}$ annihilates $e_1-e_{2n-1}$.  On the other hand, we can expand
$(e_1-e_{2n-1})P_\lambda(\, ;q,t)$ as a linear combination of Macdonald
polynomials using the Pieri identity; if we throw out those polynomials
annihilated by $L$, at most two terms remain.  Together with the identity
\[
L(P_{1^n+\lambda}(\, ;q,t))=L(e_{2n}P_\lambda(\, ;q,t))=L(P_\lambda(\, ;q,t)),
\]
we obtain an identity of the form 
\[
L(P_{\mu^2}(\, ;q,t)) = C_{\mu/\nu}(\, ;q,t) L(P_{\nu^2}(\, ;q,t)),
\]
where $\nu$ is obtained by removing a single square from the diagram of
$\mu$.  The claim then follows by induction in $|\mu|$.
\end{proof}

\begin{rem}
For many of the vanishing integrals considered below, either the linear
functional fails to factor through a homomorphism, or the homomorphism it
does factor through does not have any useful elements in its kernel, and we
thus cannot apply the Pieri identity to obtain the nonzero values.
\end{rem}

The final flaw in the above argument is that it {\em only} applies to the
symplectic case of the vanishing integral.  The point is that in the
symplectic case, the condition on compositions $\nu$ translates to a very
simple condition on the action of translations for $t=1$, namely that
certain translations should act in the same way on the monomial $x^\nu$.
For the orthogonal case, the corresponding condition on eigenvalues of
translations is actually Zariski dense; in particular, it cannot be
detected by any 
finitely generated ideal
of the Hecke algebra.  It turns out,
however, that one can deduce the orthogonal vanishing integrals from the
symplectic vanishing integral, using the fact that both can be viewed as
statements about the algebra of symmetric functions, related by a slightly
modified Macdonald involution.  (Note, in particular that 
the conjugate partition to one of the form $\mu^2$ is of
the form $2 \nu$.)

We thus obtain the following corollary; for the details, see section 8 of
\cite{bcpoly}.  Each of the four integrals is with respect to an
appropriate special case of the normalized Koornwinder density; we denote
such an $n$-dimensional integral as
$I^{(n)}_K(f;q,t;a,b,c,d)$.

\begin{cor}
For all integers $n\ge 0$ and partitions $\lambda$ with at most $n$ parts,
\[
\frac{1}{2}
I^{(n)}_K(P_\lambda(\dots,z_i^{\pm 1},\dots;q,t);q,t;\pm 1,\pm \sqrt{t})
+
\frac{1}{2}
I^{(n-1)}_K(P_\lambda(\dots,z_i^{\pm 1},\dots,\pm 1;q,t);q,t;\pm t,\pm \sqrt{t})
=
0
\]
unless $\lambda$ is of the form $2\mu$, in which case the value is
\[
\frac{C^0_{\mu}(t^{2n};q^2,t)C^-_\mu(q;q^2,t)}
     {C^0_\mu(q t^{2n-1};q^2,t)C^-_\mu(t;q^2,t)}
\]
Similarly,
\[
\frac{1}{2}
I^{(n)}_K(P_\lambda(\dots,z_i^{\pm 1},\dots,1;q,t);q,t;t,-1,\pm \sqrt{t})
+
\frac{1}{2}
I^{(n)}_K(P_\lambda(\dots,z_i^{\pm 1},\dots,-1;q,t);q,t;1,-t,\pm \sqrt{t})
=
0
\]
unless $\lambda$ is of the form $2\mu$, in which case the value is
\[
\frac{C^0_{\mu}(t^{2n+1};q^2,t)C^-_\mu(q;q^2,t)}
     {C^0_\mu(q t^{2n};q^2,t)C^-_\mu(t;q^2,t)}
.
\]
\end{cor}

\begin{rem}
Note that again the nonzero values follow via an application of the Pieri
identity from the fact that the linear functionals vanish where required.
Etingof (personal communication) has pointed out a direct proof of the
orthogonal vanishing integrals in the Jack polynomial limit, using the
construction of \cite{EK} for Jack polynomials; presumably
the Macdonald polynomial analogue  of the construction can be used to obtain
the orthogonal vanishing integral for Macdonald polynomials.  Etingof's
argument also gives no information about the nonzero values, and just as
the nature of the Hecke algebra made it impossible to use our argument in
the orthogonal case, the nature of the Etingof-Kirillov construction of
Jack polynomials makes it impossible to use Etingof's argument in the
symplectic case.
\end{rem}

\section{Other Vanishing Integrals}
\label{sec-other}

In this section, we list the remaining 
vanishing results.
For each result, we list the functional $L$ that gives the vanishing
integral, and the
associated right ideal $I$ in the Hecke algebra (of type $A$ or $C$)
that kills $L$.
We also give the subgroup $S$ of the braid group that leaves the
functional invariant in the classical limit, as this motivated
many of the relevant definitions. 
In fact, in each case the subgroup $S$ of the braid group
lies over the commutator of an involution in the extended
affine Weyl group. 
For instance,
each ideal $I$
is generated by elements
\[
U_\sigma-\chi(\sigma)
\]
where $\sigma$ is a generator of $S$ and
 $\chi$ is the character of $B(\widetilde W)$
given by its action on the constant polynomials.

In each case we argue as in section \ref{sec-sp.vanish},
that is, we exhibit a specialization of the parameters such that 
(a) in that specialization a nonzero functional  annihilated by $I$ exists
only if $\lambda$ is of the stated form, and
(b)
if $\lambda$ is of the stated form then there is unique
such nonzero functional (which can be obtained by specializing the
appropriate integral).
Since the space of functionals annihilated by an ideal can
only get bigger under specialization, this implies
(a) that if $\lambda$ is {\em not} of the appropriate form,
then for generic parameters no such functional exists, and
(b) if $\lambda$ is of the appropriate form,
such a functional exists and is generically unique. 

Note that this uniqueness is on a partition by partition basis.
It is quite possible (and indeed we give examples below)
for there to be multiple nice functionals on the space of
{\it all\/} polynomials that are all
annihilated by the same ideal
and thus satisfy  the same vanishing conditions.
(See for instance sections \ref{sec-A2}, \ref{sec-C1} below.)
For the $\widetilde S_n$ cases, this specialization is simply
$t=1$.
For the $\widetilde C_n$ cases, we must moreover take
$a=1$, $b=-1$, $c=\sqrt q$, $d=-\sqrt q$.
In all cases this has the effect of turning the Hecke algebra
into a group algebra and the (nonsymmetric) Macdonald
and Koornwinder polynomials into monomials and the density
trivial, at which point the functional is easy to evaluate.

One consequence of this global non-uniqueness is that in order
to determine the nonzero values of such a functional, it is
not enough to know how the affine Hecke algebra acts.
One must in fact consider more carefully the explicit
structure of the functional, e.g., as in the Pieri trick
used above 
(or more generally, how the double affine Hecke algebra
interacts with the functional).
Even if such a calculation could be pushed through,
this would still leave the nontrivial task of deducing
the values  on the symmetric Macdonald and Koornwinder
polynomials from the nonsymmetric ones.

In each case there is an associated family of chambers
(see section \ref{sec-comm})
such that the elements of the form $Y_\nu^C$ contained
in $S$ imply the appropriate vanishing theorem.
The fact that these chambers are not the standard chamber
implies that the standard nonsymmetric Macdonald polynomials
do not satisfy vanishing results, but the nonstandard
$E_\lambda^C$ do.
This is so because a different choice of chamber $C$ twists the Hecke algebra
module by an inner automorphism that results in an isomorphic module,
which is easily seen from the fact that
the irreducible modules $\LL_\lambda$ we consider
have distinct central characters.
In any event, while the nonstandard choice of chamber $C$
gives a different family of nonsymmetric Macdonald polynomials $E_\lambda^C$
(as they are the eigenfuctions of commuting operators $Y^C$)
the symmetric Macdonald polynomials stay the same.  ($P_\lambda^C = P_\lambda$;
symmetric functions in the $Y^C$ are always central.)
Further discussion on $E_\lambda^C$, $Y^C$ are in section \ref{sec-comm}.

\subsection{Macdonald polynomial results: $\tilde{S}_{2n}$}

\subsubsection{Case 1}
This case was discussed in section \ref{sec-sp.vanish} above.

\begin{thm}
\label{thm-vanish1.other}
\[
\label{Sn1}
\int
P_\lambda(\dots z_i^{\pm 1}\dots;q,t)
\prod_{1\le i\le n} \frac{(z_i^{\pm 2};q)}{(t z_i^{\pm 2};q)}
\prod_{1\le i<j\le n} \frac{(z_i^{\pm 1} z_j^{\pm 1};q)}
                           {(t z_i^{\pm 1} z_j^{\pm 1};q)}
dT
=
0
\]
unless $\lambda=\mu^2$, in which case (when suitably normalized) it is
\[
\frac{C^0_\mu(t^{2n};q,t^2)C^-_\mu(qt;q,t^2)}
     {C^0_\mu(qt^{2n-1};q,t^2)C^-_\mu(t^2;q,t^2)}.
\]
\end{thm}

This is the case $T=t^n$ of the symmetric function identity
\[
I_K(P_\lambda(;q,t);q,t,T;\pm\sqrt{t},\pm\sqrt{qt})
=
\begin{cases}
\frac{C^0_\mu(T^2;q,t^2)C^-_\mu(qt;q,t^2)}
     {C^0_\mu(qT^2/t;q,t^2)C^-_\mu(t^2;q,t^2)} & \lambda=\mu^2\\
0&\text{otherwise.}
\end{cases}
\]
The action of the Macdonald involution on lifted Koornwinder polynomials
dualizes this to
\[
I_K(P_\lambda(;q,t);q,t,T;\pm 1,\pm \sqrt{t})
=
\begin{cases}
\frac{C^0_\mu(T^2;q^2,t) C^-_\mu(q;q^2,t)}
{C^0_\mu(qT^2/t;q^2,t) C^-_\mu(t;q^2,t)}&\lambda=\mu^2\\
0&\text{otherwise.}
\end{cases}
\]
Taking $T\in \{t^n,t^{n+1/2}\}$ gives that the four integrals
\begin{align}
&I^{(n)}_K(P_\lambda(x_1^{\pm 1},\dots, x_n^{\pm
1};q,t);q,t;\pm 1,\pm\sqrt{t})\notag\\
&I^{(n-1)}_K(P_\lambda(x_1^{\pm 1},\dots, x_{n-1}^{\pm
1},1,-1;q,t);q,t;\pm t,\pm\sqrt{t})\notag\\
&I^{(n)}_K(P_\lambda(x_1^{\pm 1},\dots, x_n^{\pm
1},1;q,t);q,t;t,-1,\pm\sqrt{t})\notag\\
&I^{(n)}_K(P_\lambda(x_1^{\pm 1},\dots, x_n^{\pm
1},-1;q,t);q,t;1,-t,\pm\sqrt{t})\notag
\end{align}
vanish unless all ($2n$ or $2n+1$, as appropriate) parts of $\lambda$ have
the same parity.

We take 
\[
S = \langle U_0,U_n,U_i U_{2n-i}^{-1}\rangle \subseteq B(\widetilde S_{2n})
\]
The relevant chambers are such that $r$ and $r^\omega$ have the same sign, where
$\omega$ is the longest element of $S_n$, and $r$ is a root such that
$r+r^\omega\ne 0$. 
 Invariant functionals ($L(\sigma f)=0$, $\sigma\in S$)
vanish on $E^C_\lambda$ unless $\lambda_i=\lambda_{2n+1-i}$, $1\le i\le n$.

The associated right ideal of $\Hnn$ is given by
$I_1^S = \langle T_0-t,\, T_n-t \, ,T_i -  T_{2n-i}\rangle$.

The functional, which obeys $ L I_1^S  = 0$
(equivalently $L(\sigma f)=0$, $\forall  \sigma\in S$) is
\[
L(p)
=
\int p(z_1,z_2,\dots z_n,1/z_n,\dots 1/z_2,1/z_1)
\prod_{1\le i\le n} \frac{(z_i^2,q z_i^{-2};q)}{(t z_i^2,q t z_i^{-2};q)}
\prod_{1\le i<j\le n} \frac{(z_i z_j^{\pm 1},q z_i^{-1} z_j^{\pm 1};q)}
                           {(t z_i z_j^{\pm 1},qt z_i^{-1} z_j^{\pm 1};q)}
dT.
\]

\subsubsection{Case 2}
\label{sec-A2}
\begin{thm}
\label{thm-vanish2}
Let $\mu, \nu$ be partitions. Then
\[
\int
P^{(2n)}_{\mu\bar\nu}(\dots \pm \sqrt{z_i}\dots;q,t)
\prod_{1\le i\le n} \frac{((z_i/z_j)^{\pm 1};q^2)}{(t^2 (z_i/z_j)^{\pm 1};q^2)}
dT
=
\int
P^{(2n)}_{\mu\bar\nu}(\dots \pm \sqrt{z_i}\dots;q,t)
\widetilde \Delta_S^{(n)}(q^2, t^2)
dT
=
0
\]
unless $\mu=\nu$, when (suitably normalized) the integral is
\[
\frac{(-1)^{|\mu|}C^-_\mu(q;q,t) C^+_\mu(t^{2n-2} q;q,t) C^0_\mu(t^n,-t^n;q,t)}
{C^-_\mu(t;q,t) C^+_\mu(t^{2n-2} t;q,t) C^0_\mu(q t^{n-1},-q t^{n-1};q,t)}.
\]
\end{thm}
\begin{rem}
The nonzero values can be computed using the Pieri identity
as in the proof of Corollary \ref{cor-nonzeroSp}.
\end{rem}
Note that this is well-defined because
$P^{(2n)}_{\mu\bar\nu}(\dots \pm \sqrt{z_i}\dots;q,t)$
is invariant under $\sqrt{z_i} \mapsto - \sqrt{z_i}$
and it is therefore in $\C(q,t)[{z_1}^{\pm 1},\ldots, {z_n}^{\pm 1} ].$

Since multiplying a Macdonald polynomial 
by $(z_1 z_2 \cdots z_n)^m$ has the effect of addint $m$ to each part 
(which works for all $m \in \frac 1n \Z$) we can restate this in
terms of ordinary Macdonald polynomials as follows:

\begin{cor}
\[
[P_{m^n}(;q^2,t^2)] P_\lambda([2p_{k/2}];q,t)
=
0
\]
unless $\lambda=(2m)^{2n}-\lambda$.  
\end{cor}
This statement is self-dual.

We take 
\[
S = \langle U_{2i-1},U_{2i}U_{2i-1}U_{2i+1}^{-1}U_{2i}^{-1},
                     U_{2i}^{-1}U_{2i-1}U_{2i+1}^{-1}U_{2i},\shift^2\rangle.
\]
Chambers are such that $r$ and $r^\iota$ have opposite signs, where
$\iota(2i-1)=2i$, $\iota(2i)=2i-1$. 
\[
I_2^S = \langle T_{2i-1}-t \, ,T_{2i}(T_{2i-1} - T_{2i+1}),
\,\shift^2 - 1\rangle
\]
$I_2^S$-invariant functionals vanish on
$E^C_\lambda$ unless $\lambda_{2i-1}+\lambda_{2i}=0$, $1\le i\le n$.

The functional, which obeys $ L I_2^S  = 0$ is
\[
L(p)
=
\int p(z_1,-z_1,z_2,-z_2,\dots z_n,-z_n)
\prod_{1\le i<j\le n} \frac{(z_i^2/z_j^2,q^2 z_j^2/z_i^2;q^2)}
                           {(t^2 z_i^2/z_j^2,q^2 t^2 z_j^2/z_i^2;q^2)}
dT
\]

\begin{thm}
\label{thm-A2b}
Let $\lambda$ be a weight of the double cover of $GL_{2n}$,
i.e. a half-integer vector such that $\lambda_i - \lambda_j \in \Z
\, \forall i,j$.
\[
\int
P^{(2n)}_{\lambda}(\dots t^{\pm 1/2}z_i\dots;q,t)
\prod_{1\le i<j\le n} \frac{((z_i/z_j)^{\pm 1};q)}
                           {(t^2 (z_i/z_j)^{\pm 1};q)}
dT
=
0
\]
unless $\lambda_i = - \lambda_{2n+1-i}$. 
\end{thm}
We allow half-integral $\lambda$ here in order to
allow $m$ odd in the symmetric function analogue:
\begin{cor}
\[
[P_{m^n}(;q,t^2)] P_\lambda([p_k (t^{k/2}+t^{-k/2})];q,t)
=
0
\]
unless $\lambda=m^{2n}-\lambda$.
Dually,
\[
[P_{m^n}(;q^2,t)] P_\lambda(;q,t)
=
0
\]
unless $\lambda=(2m)^n-\lambda$. 
\end{cor}
\begin{rem}
Experimentally, 
the nonzero values appear to be nice, but the
kernel of the specialization 
$f \mapsto f(\cdots t^{\pm 1/2} z_i \cdots )$
is too complicated for us
to obtain recurrences from the  Pieri  identity. 
\end{rem}

For this vanishing integral, we can take the same $S_2$ and $I_2^S$,
but use a different functional:
\begin{multline}
L(p)
=
\int
p(z_1,t z_1,z_2,t z_2,\dots z_n,t z_n)
\prod_{1\le i<j\le n} \frac{(z_i/z_j,q z_j/z_i;q)}
                           {(t^2 z_i/z_j,q t^2 z_j/z_i;q)}
dT
\\
=
\int
p(z_1,t z_1,z_2,t z_2,\dots z_n,t z_n)
\widetilde \Delta_S^{(n)}(q, t^2)
dT
\end{multline}

\subsubsection{Case 3} 
\label{sec-A3}

\begin{thm}
\label{thm-vanish3}
($q\mapsto q^2$)
\[
\int
P^{(2n)}_{\mu\bar\nu}(\dots q^{\pm 1/4}z_i\dots;q,t)
\prod_{1\le i<j\le n} \frac{((z_i/z_j)^{\pm 1};q^{1/2})}
                           {(t (z_i/z_j)^{\pm 1};q^{1/2})}
dT
=
0
\]
unless $\mu=\nu$.
\end{thm}
\begin{cor}
For any partition $\lambda$,
\[
\label{eq-5.19}
[P_{m^n}(;q,t)] P_\lambda([p_k (q^{k/2}+q^{-k/2})];q^2,t)
=
0
\]
unless $\lambda=m^{2n}-\lambda$.
  Dually,
\[
[P_{m^n}(;q,t)] P_\lambda(;q,t^2)
\]
unless $\lambda=(2m)^n-\lambda$. 
\end{cor}
\begin{rem}
That \eqref{eq-5.19} holds when $\ell(\lambda) > 2n$ follows
immediately from the fact that Macdonald polynomials are triangular
with respect to the dominance order and the way the
specialization acts on monomials.

Again, the nonzero 
values appear nice, but the Pieri trick fails.
\end{rem}

We take
\[
S = \langle U_i U_{i+n}^{-1},\shift \rangle
\]
Chambers are such that $r$ and $r^\iota$ have opposite signs, where
$\iota(i)=i+n$, $\iota(i+n)=\iota(i)$.
The associated right ideal is
\[
I_3^S = \langle T_i -  T_{i+n}, \, \shift -1  \rangle.
\]
 $I_3^S$-invariant functionals vanish
on $E^C_\lambda$ unless $\lambda_i+\lambda_{i+n}=0$.

The functional is:
\begin{multline}
L(p)
=
\int
p(q^{1/2} z_1, q^{1/2} z_2,\dots q^{1/2} z_n,
  z_1,z_2,\dots z_n)
\prod_{1\le i<j\le n} \frac{(z_i/z_j,q^{1/2} z_j/z_i;q^{1/2})}
                           {(t z_i/z_j,q^{1/2} t z_j/z_i;q^{1/2})}
\\
=
\int
p(q^{1/2} z_1, q^{1/2} z_2,\dots q^{1/2} z_n,
  z_1,z_2,\dots z_n)
\widetilde \Delta_S^{(n)}(\sqrt q, t).
\end{multline}

\subsection{Koornwinder polynomial results: $\widetilde{C}_{2n}$}

\subsubsection{Case 1}
\label{sec-C1}
\begin{thm}
\label{thm-vanish1C}
In symmetric function terms,
\[
I_K(\tilde{K}_\lambda([2p_{k/2}];q,t,T;a,-a,c,-c);q^2,t^2,T;-t,-qt,a^2,c^2)
=
0
\]
unless $\lambda=\mu^2$, when it is
\[
\frac{
(-1)^{|\mu|}
C^-_\mu(qt;q,t^2)
C^+_\mu(a^2 c^2 T^2/t^4;q,t^2)
C^0_\mu(T,-a^2 T/t,-c^2 T/t,a^2 c^2 T/t^2;q,t^2)
}{
C^+_\mu(a^2 c^2 T^2/qt^3;q,t^2)
C^-_\mu(t^2;q,t^2)
C^0_{\mu^2}(a^2 c^2 T^2 q/t^2;q^2,t^2)
}.
\]
Dually,
\[
I_K(\tilde{K}_\lambda([2p_{k/2}];q,t,T;a,-a,c,-c);q^2,t^2,T;-1,-t,a^2,c^2)
=
0
\]
unless $\lambda=2\mu$, when it is
\[
\frac{
(-1)^{|\mu|}
C^-_\mu(q;q^2,t)
C^+_\mu(a^2 c^2 T^2/t^3;q^2,t)
C^0_\mu(T,-a^2 T/t,-c^2 T/t,a^2 c^2 T/t^2;q^2,t)
}{
C^-_\mu(t;q^2,t)
C^+_\mu(a^2 c^2 T^2/qt^2;q^2,t)
C^0_{2\mu}(a^2 c^2 T^2/t^3;q^2,t^2)
}.
\]
\end{thm}

The nonzero values are computed via the Pieri identities for Koornwinder
polynomials \cite{vanDiejen}.  For $T=t^{2n}$, both formal integrals become
actual integrals; similarly, for $T=t^{2n+1}$, the second formal integral
becomes:
\[
I^{(n)}_K(
K^{(2n+1)}_\lambda(\dots,\pm z_i,\dots,\sqrt{-1};q,t;a,-a,c,-c)
   ;q^2,t^2;-t,-t^2,a^2,c^2).
\]

For this identity, we work with the case $b=-a, d=-c$ of $\HnC$ and
its polynomial representation, and take
\[
S = \langle
U_{2i-1},U_{2i}^{\pm 1}U_{2i-1}U_{2i+1}^{-1}U_{2i}^{\mp 1},
U_0^{\pm 1} U_1 U_0^{\mp 1},
U_{2n}^{\pm 1} U_{2n-1} U_{2n}^{\mp 1}
\rangle
\subseteq B(\widetilde C_{2n})
\]
with associated right ideal
\[
I_1^K = \langle
T_{2i-1} -t ,T_{2i}(T_{2i-1}- T_{2i+1}),
T_0( T_1  -t),
T_{2n}( T_{2n-1}  -t )
\rangle
\]

The functional is:
\[
L(p)
=
\int
p(z_1^{1/2},-z_1^{1/2},z_2^{1/2},-z_2^{1/2},\dots z_n^{1/2},-z_n^{1/2})
\Delta_K^{(n)}(a^2, -t, c^2, -qt; q^2, t^2)
\]


\begin{thm}
In symmetric function terms,
\[
\label{IK}
I_K(\tilde{K}_\lambda([p_k(t^{k/2}+t^{-k/2})];q,t,T;a,b,c,d);
    q,t^2,T;t^{1/2}a,t^{1/2}b,t^{1/2}c,t^{1/2}d)
=
0
\]
unless $\lambda=\mu^2$.
The dual statement is:
\[
\label{IKdual}
I_K(\tilde{K}_{\lambda}(;q,t;T;a,b,c,d);q^2,t,T;a,b,c,d)
=
0
\]
unless $\lambda=2\mu$.
\end{thm}

In the case $n=1$ the above identity becomes an
identity of Askey-Wilson polynomials and admits a direct hypergeometric
proof (Rahman, personal communication).

Once again, the Pieri trick fails, but in fact the nonzero values
\[ 
\frac{
t^{-|\mu|}
C^0_\mu(T,Tab/t,Tac/t,Tad/t,Tbc/t,Tbd/t,Tcd/t,Tabcd/t^2;q,t^2)
C^+_\mu(T^2abcd/t^4;q,t^2)
C^-_\mu(qt;q,t^2)
}{
C^0_{2\mu^2}(T^2abcd/t^2;q,t^2)
C^+_\mu(T^2abcd/qt^3;q,t^2)
C^-_\mu(t^2;q,t^2) 
}
\]
for the first integral and
\[
\frac{
q^{|\mu|}
C^0_\mu(T,Tab/t,Tac/t,Tad/t,Tbc/t,Tbd/t,Tcd/t,Tabcd/t^2;q^2,t)
C^+_\mu(T^2abcd/t^3;q^2,t)
C^-_\mu(q;q^2,t)
}{
C^0_{2\mu^2}(T^2abcd/t^2;q^2,t)
C^+_\mu(T^2 abcd/t^2q;q^2,t)
C^-_\mu(t;q^2,t)
}.
\]
for the second can be obtained as a limit of the elliptic version derived
in \cite{littlewood}.

We take $S$ and $I_1^K$ as above,  but now with 
generic $a,b,c,d$ and  the functional we need is:
\[
L(p)
=
\int
p(t^{-1/2} z_1,t^{1/2} z_1,t^{-1/2} z_2,t^{1/2} z_2,\dots
t^{-1/2} z_n,t^{1/2} z_n)
\Delta_K^{(n)}(t^{1/2}a, t^{1/2}b, t^{1/2}c, t^{1/2}d; q, t^2)
\]

\subsubsection{Case 2}  
\label{sec-C2}
\begin{thm}
\label{thm-vanish2C}
In symmetric function terms,
\[
\label{IK2}
I_K(\tilde{K}_\lambda([p_k(q^{k/2}+q^{-k/2})];q^2,t,T^2;a,b,qa,qb)
   ;q,t,T;\pm\sqrt{t},q^{1/2}a,q^{1/2}b)
=
0
\]
unless $\lambda=\mu^2$.

The dual statement is:
\[
I_K(\tilde{K}_\lambda( \,;q,t^2;T^2;a,b,ta,tb)
   ;q,t,T;\pm \sqrt{t},a,b)
=
0
\]
unless $\lambda=2\mu$.
\end{thm}
We take 
$c=q^{1/2}a$, $d=q^{1/2}b$ (so consider the case $q \mapsto \sqrt q$
above)
\[
S = \langle
U_0U_{2n}^{-1},U_iU_{2n-i}^{-1},U_n\rangle
\subseteq B(\widetilde C_{2n})
\]
and associated right ideal
\[
S = \langle
T_0 - T_{2n} \, ,T_i - T_{2n-i},T_n -t \rangle
\]

The functional is:
\[
L(p) :=
\int
p(q^{1/4} z_1,\dots q^{1/4} z_n,q^{1/4}/z_n,\dots q^{1/4}/z_1)
\Delta_K^{(n)}(\sqrt{t}, -\sqrt{t}, q^{1/4}a, q^{1/4}b; q^{1/2}, t).
\]



\section{A construction using the Hecke algebra}
 
\label{sec-hecke}

In this section, we give another proof of the existence of
nonzero functional $L$ in the non-vanishing case 
(another proof of the vanishing condition can also be deduced)
 along with an
 explicit construction of $L \in \Llam$ (up to scalar). 
We only do this for the $\widetilde S_{2n}$ case, leaving
the Koornwinder case to the reader.
We do not explicitly compute the scalar that relates the $L$ constructed
in this section to the integral given in section \ref{sec-sp.vanish}.
We also warn the reader that since we are computing in 
$(\Vquot)^*$ versus $\Vlam^*$, we do not give information about
$L(E_\mu)/L(E_\nu)$ except when $\mu^+ = \nu^+$.

In what follows we will use the presentation of $\Hnn$ as generated by
$T_1, T_2, \ldots, T_{2n-1}$, $Y_1^{\pm 1}, \ldots, Y_{2n}^{\pm 1}$
because it allows us to work more explicitly with a basis of 
$\Llam$ given by simultaneous $Y_i$-eigenfunctionals.

This presentation also gives us another description of
$\Vquot$ and of $\Llam = (\Vquot)^*$. 
Given $\lambda \vdash 2n$, let 
$J =J_\lambda=  \{ j \mid s_j \lambda = \lambda \}$,
let $\Hlam$ be the parabolic subalgebra generated by 
$\{ T_j \mid j \in J \}$ and all the 
$Y_1^{\pm 1}, \ldots, Y_{2n}^{\pm 1}$, and
let $\Clam$ be the
one-dimensional $\Hlam$ module on which $Y_i - q^{\lambda_i} t^{2n -i} =0$,
$T_j -t =0 \, \forall j \in J$.
Then we have
$$\Llam \isom  \Clam \otimes_{\Hlam} \Hj{2n}.$$
(Note that $\Vquot$ is
isomorphic to $\Hj{2n} \otimes_{H( w_0 \lambda)} \C(q,t)_{
w_0 \lambda}$ 
(which is isomorphic to 
$\Hj{2n} \otimes_{H(  \lambda)} \C(q,t)_{  \lambda}$ 
when $q,t$ are generic)
and thus  is  in this sense self dual.
This can be seen directly or follows from the Mackey decomposition
of $V$.)

For ease of notation, we introduce the standard invariant
form $\enbrac{\,\,}{\,}$,
and
let $\delta = \delta_{2n} = (2n-1, \ldots, 2, 1, 0)$, $\epsi =
(0, \ldots, 0, 1, 0, \ldots, 0), \alpha_i = \epsi - \eps{i+1}$.
We can then write $\lambda_i = \enbrac{\lambda}{\epsi}$.
We also have 
$\enbrac{  \lambda}{\mu} = \enbrac{w\lambda}{w \mu}$,
where $w \in S_{2n}$ acts as $w (\lambda_1, \ldots, \lambda_{2n}) =
(\lambda_{w^{-1}(1)}, \ldots \lambda_{w^{-1}(2n)}).$

We observe that the center $Z(\Hnn)$ is
given by symmetric Laurent polynomials in $ Y_1, \ldots, Y_{2n}$
and each $\Llam$ has distinct central character.
Further, the $Y$-weight spaces of $\Llam$ are all one-dimensional
and hence give a distinguished basis of the module, up to scalars.
From the above  description of $\Llam$,
it is easy to see
that basis of
simultaneous $Y_i$-eigenvectors is $\{ \vw \mid w \in \WJ\}$ with
$$ \vw (Y_i - q^{\enbrac{w^{-1}\lambda}{\epsi}}
t^{\enbrac{w^{-1}\delta}{\epsi}}) = 0,$$
where $\WJ$ is the set of minimal length right coset representatives
for $\langle s_j \mid j \in J \rangle \subseteq S_{2n}$.
We normalize this basis so that the {\it right\/}
action of the $T_i, 1 \le i < n$ is given by 
\[ 
\vw T_i = \frac{t-1}{1- q^{- \enbrac{w^{-1} \lambda}{\alpha_i} }
t^{- \enbrac{w^{-1} \delta}{\alpha_i} } } \vw 
+ 
 \frac
{ 1- q^{ \enbrac{w^{-1} \lambda}{\alpha_i} }
     t^{1+\enbrac{w^{-1} \delta}{\alpha_i} } } 
{1- q^{- \enbrac{w^{-1} \lambda}{\alpha_i} }
    t^{- \enbrac{w^{-1} \delta}{\alpha_i}  } } \vwsi 
\]
with the convention $\vwsi = 0$ if $w s_i \not\in \WJ$. 
In that case, notice $\vw T_i = t \vw$, $\vw Y_iY_{i+1}^{-1} = t \vw$,
and in particular $\lambda_{w(i)} = \lambda_{w(i+1)}$.
Observe the above action does not depend on the relative lengths
$\ell(w)$ and $\ell(w s_i)$, which is why this particular normalization
is preferred in this setting. 

We note that this basis
dual to the one given by
the nonsymmetric Macdonald polynomials 
up to proportionality. 
(We leave it to the reader to rescale as necessary, possibly also
rescaling the $T_i$, to get exactly the dual basis.)

We want to express our functional $L$ (given by integrating a specialized
polynomial against a given density) in terms of this basis
$\{ \vw \}$.
To do this, we conjugate the right ideal $I$ such that $L \cdot I = 0$
to a related right ideal $T_u I T_u^{-1}$ which has 
a nicer presentation in terms of the $Y_i$.
This corresponds to working with the twisted module  $\Llam^u \isom \Llam$.
Hence we explicitly describe the functional $L T_u^{-1}$ in terms
of the dual basis to nonsymmetric Macdonald polynomials, but not $L$ 
itself. 
The vanishing result stated above will hold for suitable 
nonsymmetric Macdonald polynomials with $E_\nu$ replaced
by $T_u E_\nu$. 
This twist by $u$ is motivated by the viewpoint of section
\ref{sec-comm} regarding nonstandard large commutative
subalgebras.

For each functional $L$ and corresponding ideal $I$ such that $LI=0$,
we describe $T_u I T_u^{-1} = I'$, determine all $\lambda$ such that
there exists nonzero $v \in \Llam$ with $vI'=0$, and show this
$v$ is unique up to scalar.


It will  follow from our explicitly computed generators
of $T_u I T_u^{-1}$ that it contains a large binomial ideal in
the commutative subalgebra
$\C[Y_i^{\pm 1}]$.
This translates directly to conditions under which $L T_u^{-1}( E_\nu)$
is forced to vanish.
In particular this implies for any $L$ such that $L T_u^{-1} I = 0$
that $L T_u^{-1} (E_\nu)$ vanishes.
In each case this will immediately give the desired vanishing
result for $P_{\nu^+}$. 
However, one can ask for something stronger, namely that
for each partition $\lambda$ either the stated vanishing
condition holds or there exists a {\it unique \/} $I$-killed
functional. 

In what follows, all ideals are right ideals. 

\subsection{Second proof of Theorem \ref{thm-vanish1.other} }
\label{sec-hecke-ideal1} 
Recall $I_1^S = \< T_0-t, \;  T_n -t,\;  T_i - T_{2n-i} \, (1\le i < n) \>$
is the right ideal with given generators. 
Let $u$ be the permutation defined by 
$$ u(i)=\begin{cases} 2(n-i) +1 \quad  & i\le n \\
		2(i-n) \quad & i > n \end{cases}$$
and set $I'= {I_1'}^{A} = T_u I_1^S T_u^{-1}$.
Observe $\ell(u)  = 2 \binom{n}{2}$.
Then
$$
I' = \<  T_{2i-1} -t,  \; T_{2i}(T_{2i+1} - T_{2i-1}), \;  t Y_{2i} - Y_{2i-1}
\, (1\le i < n) \>, $$

One can verify
\begin{align*}
 T_u (T_i - T_{2n-i})  T_u^{-1} & = \;
T_{u(2n-i)}(T_{u(i)} - T_{u(i)-2})  T_{u(2n-i)}^{-1}
	= T_{2(n-i)}(T_{2(n-i)+1} - T_{2(n-i)-1}) T_{2(n-i)}^{-1}, \\
 T_u (T_n - t)  T_u^{-1} \quad  & = \; T_1 -t,\\
 T_u (t T_0^{-1} - 1)  T_u^{-1}& = \; Y_{2n-1}^{-1}  Y_{2n}T_{2n-1} -1.
\end{align*}
From the first two equations, we can show $T_{2i-1} -t \in I'$,
inductively as 
$T_{2i+1} -t = ( (T_{2i-1} -t)T_{2i} T_{2i-1} 
+ T_{2i}(T_{2i+1}-T_{2i-1})T_{2i}^{-1} (T_{2i}^2 - t T_{2i}) )
T_{2i+1}^{-1} T_{2i}^{-1} \in I'$.
Then $ t Y_{2n-1}^{-1}  Y_{2n} -1 =
  Y_{2n-1}^{-1}  Y_{2n} t - \bT_{2n-1} + T_{2n-1}-t
=
 ( Y_{2n-1}^{-1}  Y_{2n}T_{2n-1} -1) \bT_{2n-1} +  (T_{2n-1}-t)
\in I'.
$

To show $ t Y_{2i} - Y_{2i-1} \in I'$,
it suffices to show $Y_{2i}(T_{2i-1} -t) \in I'$
as  $Y_{2i}(T_{2i-1} -t) = (T_{2i-1} -t) Y_{2i-1} - (t Y_{2i} - Y_{2i-1}).  $
Note
\begin{align*}
Y_{2i}(T_{2i-1} -t) 
&= t^{-2} T_{2i} T_{2i+1} Y_{2i+2}  T_{2i+1} T_{2i} (T_{2i-1} -t)
\\
&= t^{-2} T_{2i} T_{2i-1} Y_{2i+2}  T_{2i+1} T_{2i} (T_{2i-1} -t)
 + t^{-2} (T_{2i} T_{2i+1} -  T_{2i} T_{2i-1}) 
	Y_{2i+2}  T_{2i+1} T_{2i} (T_{2i-1} -t)
\\
&\in t^{-2}  Y_{2i+2} T_{2i} T_{2i-1}  T_{2i+1} T_{2i} (T_{2i-1} -t)
 +  I'
\\
&= t^{-2}  Y_{2i+2} (T_{2i+1} -t) T_{2i} T_{2i-1}  T_{2i+1} T_{2i}
 +  I'
\end{align*}

Because $ t Y_{2i} - Y_{2i-1} \in I'$ 
it follows that if  a functional $L'$ is annihilated by $I'$
then $L'(E_\mu)=0$ unless $\mu_1 = \mu_2, \; \mu_3=\mu_4$ and so on,
thus directly proving the vanishing result, Theorem \ref{thm-vanish1}.

We may thus restrict our attention to partitions of the form
$\lambda = \mu^2$.
We wish to show that in this case, $\Llam$ contains a unique
$I'$-killed functional and give an explicit expression for
that functional in terms of the basis $\{\vw\}$.
Of course, it in only possible to determine the functional up
to an overall scalar (and in fact because we are only considering
this one partition at a time, we have such a scalar for every valid
partition).
What this does determine is the relative values of an $I'$-killed
functional on nonsymmetric Macdonald polynomials.
The actual values of such a functional are at least in 
principle determined by its values on symmetric Macdonald polynomials
(since for $t=1$ we can exhibit a functional for which
those values are nonzero).
Moreover, experimentally, the resulting scale factors are still nice.
However, it appears somewhat nontrivial to prove a closed form.

Next we will determine under what conditions $\Llam$ contains a functional
annihilated by $I'$, and show that it is unique up to scalar.
We will give an explicit expression for this functional in terms
of the $\vw$.

We will need some more notation.

For $w \in W$ let $R(w) = \{ \alpha > 0 \mid w \alpha < 0 \}$.
Notice for $w \in S_{2n}$ we have $R(w) = \{ \epsi - \epsj \mid i < j,
w(i) > w(j) \}$, and $|R(w)| = \ell(w)$.
For $\iota$ an involution acting on the weight lattice,
let $\Ri (w)= \{ \frac 12 (\alpha + \iota(\alpha))
\mid \alpha \in R(w) \}$. 
Since $\iota$ is an involution, the sizes of its orbits are either
one or two.  When it is necessary to differentiate, we set
$\Ri_1 = \{ \alpha \in R(w) \mid \iota(\alpha) = \alpha \}$,
$\Ri_2 = \{ \frac 12 (\alpha + \iota(\alpha)) \mid \alpha \in R(w) ,
 \iota(\alpha) \neq \alpha \}$.


\begin{prop}
Suppose $\lambda = \mu^2$.
Then any $v \in \Llam$ with $vI'=0$ is proportional to the nonzero
$I'$-killed functional
\begin{gather}
\label{eq-ideal1v}
  \mathop{\sum_{w \in \WJ }}_{\iota(w^{-1} \lambda )= w^{-1} \lambda} 
\bigl(
\prod_{\beta \in \Ri(w^{-1}) } 
  q^{-\enbrac{ \lambda}{\beta} } t^{-\enbrac{\delta}{\beta} }  
\frac{ 1- q^{ -\enbrac{ \lambda}{\beta} } t^{1-\enbrac{\delta}{\beta} } } 
     { 1- q^{\enbrac{ \lambda}{\beta} } t^{1+\enbrac{\delta}{\beta} } } 
\bigr)
\vw,
\end{gather}
where $\iota$ is the involution on the weight lattice with
$\iota(\eps{2i-1}) = \eps{2i}$.
\end{prop}

\begin{proof}
Write $v = \sum_{w \in \WJ} c_w \vw$ and suppose $v I' = 0$.

That $v ( t Y_{2i} - Y_{2i-1}) = 0$ forces $\enbrac{ w^{-1} \lambda}{\alpha_i}
=0$  and $\enbrac{ w^{-1} \delta}{\alpha_i}-1 =0$
whenever $c_w \neq 0$.
In particular $v \neq 0$ implies $\iota \lambda = \lambda$, which we have
already included in our hypotheses as $\lambda = \mu^2$.
Also, it automatically follows for such an expression
that $v(T_{2i-1} - t) = 0$.

That $v T_{2i}(T_{2i+1}-T_{2i-1}) = 0$ forces a relation on $c_w$
and $c_{w s_{2i} s_{2i+1} s_{2i-1} s_{2i} }$, and the resulting
relation between nonzero $c_w$ and $c_{id}$ is independent of 
reduced expression for $w$ and given by \eqref{eq-ideal1v}.


\end{proof}

\subsection{Second proof of Theorems \ref{thm-vanish2}, \ref{thm-A2b} }
 In order to accomodate Theorem \ref{thm-A2b},
we must allow half-integral weights, i.e., include
$(z_1 z_2 \cdots z_{2n})^{-1/2}$ in the algebra of polynomials
on which we act. 
Recall $I_2^S =
\< \shift^2 -1, \;  T_0(T_{2n-1} - T_1), \;  T_{2i-1} -t,  T_{2i}(T_{2i+1} - T_{2i-1}), \,
(1\le i < n) \>$. 
Let $v$ be the permutation defined by 
\begin{gather*}
v(2i+1) = n-i \\
v(2i) = n+i
\end{gather*}
and set $I' = T_{v} I T_{v}^{-1}$. 
Notice that $v = u^{-1}$ with $u$ the permutation for the ideal
in section \ref{sec-hecke-ideal1}, so that
$I' = T_{u^{-1}} I T_{u^{-1}}^{-1}$.
Then
$$
I' = \< T_n -t,  \;T_i - T_{2n-i}, \;    Y_{i}  Y_{2n-i+1}-t^{2n-1},
 (1\le i < n) \>, $$

We can use the same computations as with the first ideal,
using the fact there is an anti-involution $*$ on the Hecke algebra
sending $T_w \mapsto T_{w^{-1}}$, i.e., if
$T_u a T_u^{-1} = b$ then
$T_v b^* T_v^{-1} = a^*$. 
Hence we get
$ T_n -t, T_i - T_{2n-i} \in I'$.
To be more precise,
\begin{align*}
 T_i - T_{2n-i}   & = \;
T_v T_{u(2n-i)}(T_{u(i)} - T_{u(i)-2})  T_{u(2n-i)}^{-1} T_v^{-1}
=
T_v T_{2(n-i)}(T_{2(n-i) +1} - T_{2(n-i) -1})  T_{2(n-i)}^{-1} T_v^{-1}
 \\
 T_n - t\quad  & = \; T_v (T_1 -t)  T_v^{-1}.
\end{align*}

One can verify  
\begin{align*}
T_{u^{-1}} (\shift^{-2} -1) T_{u^{-1}}^{-1}
&= T_{n-1} \cdots T_2 T_1 T_n T_{n-1} \cdots T_2 \shift^{-2} T_{2n-2}^{-1}
	\cdots T_1^{-1} -1 \\
&\in
 T_{2n-1} \cdots  T_{n+1} T_n T_{n-1} \cdots T_2 \shift^{-2} T_{2n-2}^{-1}
	\cdots T_1^{-1} -1 \; + I' \\  
&= 
 t^{2n-1}Y_{2n}^{-1}  Y_{1}^{-1}-1 + I'.
\end{align*}
The second 
step comes from the fact
that $T_{2n-1} \cdots  T_{n+1}  - T_{n-1} \cdots  T_1
= \sum_{i=1}^{n-1} (T_{2n-i} - T_i ) T_{2n-i-1} T_{2n-i-2} \cdots T_{n+1}
	T_{i-1}$   $ \cdots T_2 T_1 \in I'$.

\begin{prop}
Suppose $\lambda $ satisfies
$\lambda_{i} =- \lambda_{2n - i+1 }, 
(1 \le i < n)$,
 i.e. if we set $\nu = ({\frac m2}^{2n}) + \lambda$,
then $\nu = (m^{2n}) - \nu$
(even for $m=0$).
Then any $v \in \Llam$ with $vI'=0$ is proportional to the nonzero
$I'$-killed functional
\begin{gather}
\label{eq-ideal2v}
 v = \mathop{\sum_{w \in \WJ }}_{\iota(w^{-1} \lambda) = w^{-1} \lambda} 
\bigl(
\prod_{\beta \in \Ri_2(w^{-1}) } 
(-1)  q^{-\enbrac{ \lambda}{\beta} } t^{-\enbrac{\delta}{\beta} }  
\frac{ 1- q^{-\enbrac{ \lambda}{\beta} } t^{1-\enbrac{\delta}{\beta} } } 
     { 1- q^{\enbrac{ \lambda}{\beta} } t^{1+\enbrac{\delta}{\beta} } } 
\prod_{\beta \in \Ri_1(w^{-1}) } 
(-1)  q^{-\enbrac{ \lambda}{\beta} } t^{-\enbrac{\delta}{\beta} }  
\bigr)
\vw.
\end{gather}
Here $\iota$ is the involution on the weight lattice with
$\iota(\eps{i}) = -\eps{2n-i+1}$.
\end{prop}
\begin{proof}
In the above expression $\Ri_1(w^{-1})$ represents the $\iota$-orbit sums 
on $R(w^{-1})$ where the orbit has size $1$,
and
$\Ri_2(w^{-1})$ represents the $\iota$-orbit sums 
on $R(w^{-1})$ where the orbit has size $2$.

Write $v = \sum_{w \in \WJ} d_w \vw$ and suppose $v I' = 0$.

That $v (    Y_{i}  Y_{2n-i+1}-t^{2n-1}) = 0$ forces
$\enbrac{ w^{-1} \lambda}{\epsi + \eps{2n-1}} =0$,
 $\enbrac{ w^{-1} \delta}{\epsi + \eps{2n-1}}- 2n +1 =0$
whenever $d_w \neq 0$.
In particular $v \neq 0$ implies $\iota \lambda = \lambda$, which
is in the hypotheses of our proposition.

That $v (T_{i}-T_{2n -i}) = 0$ forces a relation on $d_w$
and $d_{w s_{i} s_{2n-i}  }$ corresponding to the first term
in the above product. (Note that $\iota (\alpha_i) = \alpha_{2n-i}$.)
That $v (T_{n}-t) = 0$ forces a relation on $d_w$
and $d_{w s_{n}  }$ corresponding to the second term
in the above product. (Note that $\iota (\alpha_n) = \alpha_{n}$.)
The resulting
relation between nonzero $d_w$ and $d_{id}$ is independent of 
reduced expression for $w$ and given by \eqref{eq-ideal2v}.


\end{proof}
\subsection{Second proof of Theorem \ref{thm-vanish3} }
Again we must allow half-integral weights.

Recall $I = 
\< \shift -1,  \;   T_{i} - T_{i+n},\; 
(0\le i < n) \>$. 
Let $u$ be the permutation defined by 
$$ u(i)=\begin{cases} i \quad  & i\le n \\
		3n-i+1 \quad & i > n \end{cases}$$
and set $I' = T_u I T_u^{-1}$.
Note that $u$ is the longest element of $\underbrace{S_1
\times \cdots \times S_1}_n \times S_n$.

Then
$$
I' = \< 
T_n + 1 -t - t^{1-n} Y_{n+1},\; 
T_{i} -T_{2n-i},\;      Y_{i}  Y_{2n-i+1} - t^{2n-1},
 (1\le i < n)
\>. $$

One can verify
\begin{align*}
T_u (T_i - T_{i+n})  T_u^{-1} &= T_i-T_{2n-i}\\
T_u (\shift -1)  T_u^{-1} &= 
T_{n+1} T_{n+2} \cdots  T_{2n-1} \shift  T_{n+1}^{-1} \cdots T_{2n-1}^{-1} -1.
\end{align*}
Hence $I' \ni T_{n+1} T_{n+2} \cdots  T_{2n-1} \shift - T_{2n-1}\cdots  T_{n+1}
\equiv   T_{n+1} T_{n+2} \cdots  T_{2n-1} \shift - T_{n-1}\cdots  T_{1}$.
And so
$I' \ni T_{n+1}  \cdots  T_{2n-1} \shift T_{1}^{-1} $
$\cdots T_{n-1}^{-1} T_{n}^{-1} -T_{n}^{-1} 
= t^{-n} Y_{n+1} - T_{n}^{-1} = 
 t^{-1}( t^{1-n} Y_{n+1} -( T_{n}+1-t) )$.

Then also $I' \ni   (t^{-n} Y_{n+1} - T_{n}^{-1})(t^{n} Y_{n} + t^{2n-1}T_{n})
= 	Y_{n+1} Y_n  - t^{2n-1}.$
From this it is easy to show $Y_i Y_{2n-i+1} - t^{2n-1} \in I'.$


\begin{prop}
Suppose $\lambda $ satisfies
$\lambda_{i} =- \lambda_{2n - i+1 }, 
(1 \le i < n)$,
 i.e. if we set $\nu = ({\frac m2}^{2n}) + \lambda$,
then $\nu = (m^{2n}) - \nu$
(even for $m=0$, sorting parts in the latter expression so it is a partition).
Then any $v \in \Llam$ with $vI'=0$ is proportional to the nonzero
$I'$-killed functional
\begin{gather}
\label{eq-ideal3v}
 v = \mathop{\sum_{w \in \WJ }}_{\iota(w^{-1} \lambda) =w^{-1} \lambda} 
\bigl(
\prod_{\beta \in \Ri_2(w^{-1}) } 
(-1)  q^{-\enbrac{ \lambda}{\beta} } t^{-\enbrac{\delta}{\beta} }  
\frac{ 1- q^{-\enbrac{ \lambda}{\beta} } t^{1-\enbrac{\delta}{\beta} } } 
     { 1- q^{\enbrac{ \lambda}{\beta} } t^{1+\enbrac{\delta}{\beta} } } 
\prod_{\beta \in \Ri_1(w^{-1}) } 
(-1)  
\frac{ 1- q^{-\frac {\enbrac{ \lambda}{\beta}} 2}
	t^{\frac  {1-\enbrac{\delta}{\beta}} 2}} 
     { 1- q^{\frac  {\enbrac{ \lambda}{\beta} }2}
	t^{\frac  {1+\enbrac{\delta}{\beta}} 2} } 
\bigr)
\vw
\end{gather}
Here again $\iota$ is the involution on the weight lattice with
$\iota(\eps{i}) = -\eps{2n-i+1}$.
\end{prop}
\begin{proof}
In the above expression $\Ri_1(w^{-1})$ represents the $\iota$-orbit sums 
on $R(w^{-1})$ where the orbit has size $1$,
and
$\Ri_2(w^{-1})$ represents the $\iota$-orbit sums 
on $R(w^{-1})$ where the orbit has size $2$.

Write $v = \sum_{w \in \WJ} b_w \vw$ and suppose $v I' = 0$.

That $v (    Y_{i}  Y_{2n-i+1}-t^{2n-1}) = 0$ forces
$\enbrac{ w^{-1} \lambda}{\epsi + \eps{2n-1}} =0$,
 $\enbrac{ w^{-1} \delta}{\epsi + \eps{2n-1}}- 2n +1 =0$
whenever $b_w \neq 0$.
In particular $v \neq 0$ implies $\iota \lambda = \lambda$,
which we have already included in our hypotheses.

That $v (T_{i}-T_{2n -i}) = 0$ forces a relation on $b_w$
and $b_{w s_{i} s_{2n-i}  }$ corresponding to the first term
in the above product. (Note that $\iota (\alpha_i) = \alpha_{2n-i}$.)
That $v (T_{n}+1-t - t^{1-n} Y_{n+1}) = 0$ forces a relation on $b_w$
and $b_{w s_{n}  }$ corresponding to the second term
in the above product. (Note that $\iota (\alpha_n) = \alpha_{n}$.)
The resulting
relation between nonzero $b_w$ and $b_{id}$ is independent of 
reduced expression for $w$ and given by \eqref{eq-ideal3v}.


\end{proof}

\section{Extended affine Weyl groups}
\label{sec-extended}

Let $W$ be a finite Weyl group acting on a Euclidean space $\R^n$, with
associated root lattice $\Lambda_0$, not assumed to span $\R^n$.  A {\em
  generalized weight lattice} for $W$ is a lattice $\Lambda$ (spanning
$\R^n$) containing $\Lambda_0$ such that
\[
\frac{2\langle r,\nu \rangle}{\langle r,r\rangle}\in \Z
\]
for all roots $r$ and vectors $\nu \in\Lambda$.  (Note that there is a
one-to-one correspondence between isomorphism classes of pairs
$(W,\Lambda)$ and isomorphism classes of connected compact Lie groups; here
$\Lambda$ is the inverse image of the identity element under the
exponential map.)

An  {\em  extended  affine  Weyl  group}  is  then  a  group  of  the  form
$\widetilde{W}=G\ltimes \Lambda$,  where $\Lambda$ is  a generalized weight
lattice for a finite Weyl  group $W$, and $G\subset \Aut(\Lambda)$ contains
$W$ as a normal subgroup.
Given $\nu \in \Lambda$, we denote the corresponding element of
$\widetilde{W}$ by $\tau_\nu $ to avoid confusion.

An alcove is the closure of a fundamental region for the normal
subgroup $W\ltimes \Lambda_0$; the {\em standard alcove} is the
unique alcove containing the origin contained in the fundamental
chamber of W.  The union of the boundaries of the alcovies is a
union of hyperplanes; the distance between two alcoves is the number
of such hyperplanes that separate their interiors.  Given $w\in
\widetilde{W}$, the length of $w$ is the distance between the standard
alcove and its image under $w$.  In particular, the elements of length 0
are those that preserve the standard alcove, and there is a natural map
from $\widetilde{W}$ to the length 0 subgroup with kernel $W\ltimes\Lambda_0$.

The braid group $B(\widetilde{W})$ is generated by elements $\U(w)$ for $w\in
\widetilde{W}$, subject to the relations $\U(w_1w_2)=\U(w_1)\U(w_2)$ whenever
$\ell(w_1w_2)=\ell(w_1)+\ell(w_2)$; thus $B(\widetilde{W})$ contains a subgroup
identified with the length 0 subgroup of $\widetilde{W}$, and is generated over
this subgroup by $\U(s)$ for $s$ of length 1.

  The Hecke algebra
$H(\widetilde{W})$ is obtained from the group algebra of $B(\widetilde{W})$ by
adding further quadratic relations
$(\U(s)-t_s^{\frac 12})(\U(s)+t_s^{-\frac 12})=0$.
 We require $t_s = t_{s'}$ if
$s$ and $s'$ are conjugate, since then 
$U(s)$ and $U(s')$ are conjugate.
More generally, if $\sigma \in \widetilde{W}$
is of length $0$ 
(and so acts on the affine Dynkin diagram)
we have $\sigma$ act  on scalars  by
$$\sigma t_s \sigma^{-1} = t_{\sigma s\sigma^{-1} } .$$
If there are simple reflections which are conjugate in
$\widetilde{W}$ but not in $W\ltimes \Lambda_0$
then this action on scalars is nontrivial, and therefore
the resulting extended affine Hecke algebra is  no
longer a central algebra over $\C[t_s^{\pm 1/2}]$.
However if we specialize the $t_s$ appropriately, one can indeed
obtain 
a central algebra over $\C[t_s^{\pm 1/2}]$.
Alternatively, we can view such $\sigma$ as giving an intertwining
map between two different Hecke algebras.

For instance, in the case of $\HnC$ the outer involution $\sigma$
in general gives an intertwining map between two different instances
of $\HnC$. In particular it takes nonsymmetric Koornwinder polynomials
for one set of parameters to  nonsymmetric Koornwinder polynomials
with modified parameters.
This becomes significant because the construction of $Y$ operators
given in the next section includes such intertwiners and
this explains for instance the difference operator of \cite{bcpoly}.

For the cases $\widetilde{S}_n,$ $\widetilde{C}_n$
which are of particular interest to us, we can represent elements
of the corresponding braid groups pictorially as periodic braids.
We follow
(American) book-spine conventions; that is, the
leftmost symbol in a word corresponds to the topmost move in the
corresponding braid picture.
 To save space,
commuting symbols may be drawn as occuring at the same time.

 The
generators of the braid group are denoted $U_i$; in the Hecke algebra, they
satisfy $U_i-U_i^{-1} = t_i^{1/2}-t_i^{-1/2}$, and we define
$T_i=\sqrt{t_i} U_i$.

In $\widetilde{S}_n,$ 
$U_i $ corresponds to a picture in which
(reading down) the $j$th strand
(from the left) crosses under the $j+1$st
strand for all $j \equiv i \mod n$.
\begin{figure}[htb]
\begin{center}
\epsfig{figure=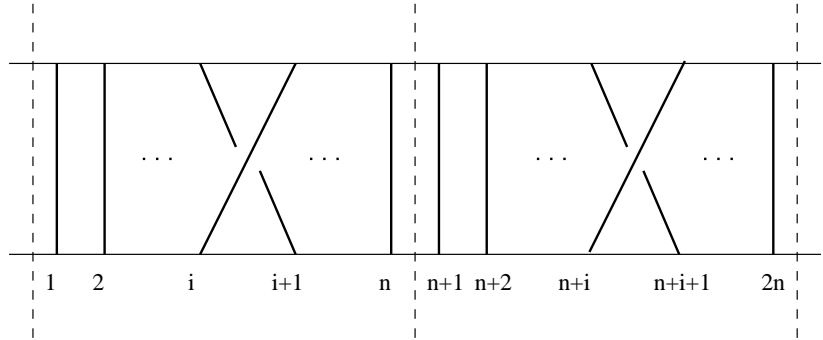}
\end{center}
\caption{$U_i \in B(\widetilde{S}_n)$ }
\end{figure}
\begin{figure}[htb]
\begin{center}
\epsfig{figure=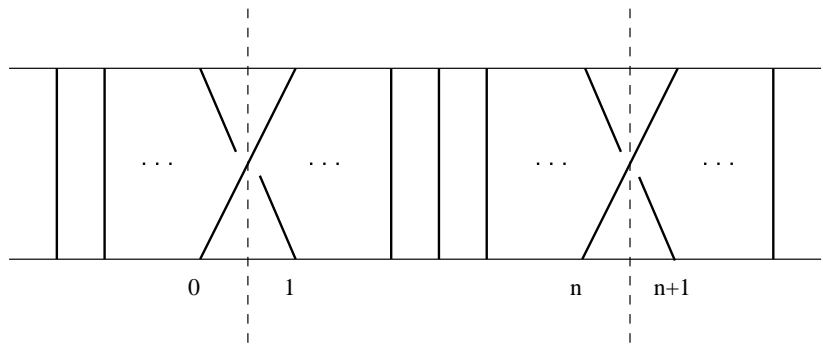}
\end{center}
\caption{$U_0 \in B(\widetilde{S}_n)$ }
\end{figure}
Similarly $\shift$ corresponds to the operation that simply
moves each strand one step to the right.
\begin{figure}[htb]
\begin{center}
\epsfig{figure=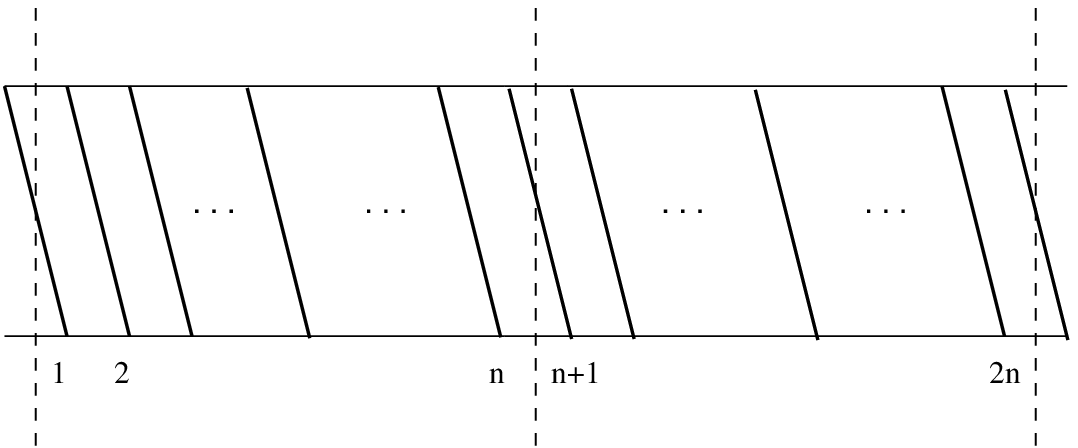}
\end{center}
\caption{$\shift \in B(\widetilde{S}_n)$ }
\end{figure}
%

\begin{figure}[htb]
\label{fig-Y1}
\begin{center}
\epsfig{figure=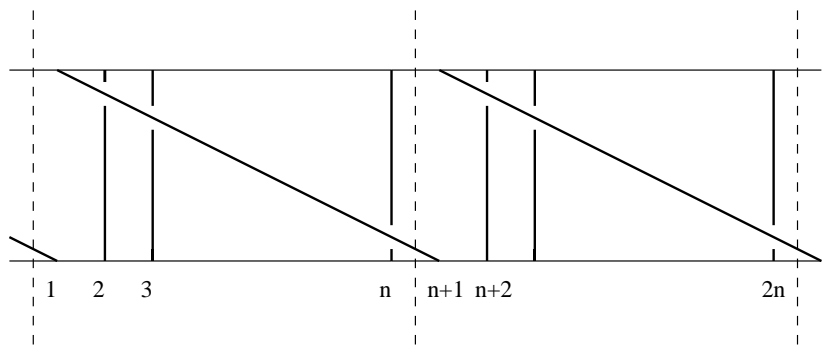}
\end{center}
\caption{$Y_1 \in B(\widetilde{S}_n)$ }
\end{figure}

\begin{figure}[htb]
\label{fig-Yi}
\begin{center}
\epsfig{figure=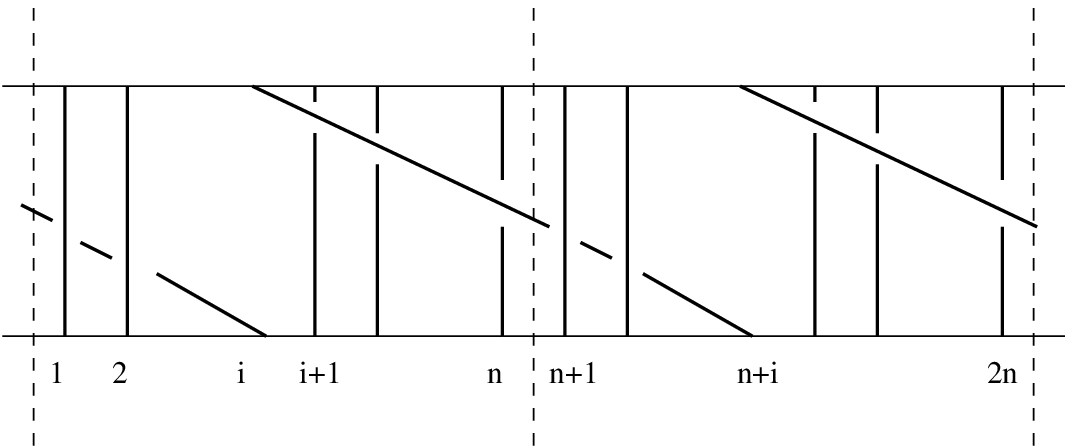}
\end{center}
\caption{$Y_i \in B(\widetilde{S}_n)$ }
\end{figure}

The elements $Y_i \in B(\widetilde{S}_n)$ (i.e., the
elements of the braid group given by replacing $T_i$ by $U_i$
and $\bT_i$ by $U_i^{-1}$ 
in \eqref{eq-Y1}, \eqref{eq-Y2}, \eqref{eq-Yn})
 moves the strands congruent to $i$ $\mod n$ $n$ steps to the right, underneath
the adjacent strands congruent to $1\dots i-1$ and over the remaining strands.
See figures \eqref{fig-Y1}, \eqref{fig-Yi}.
%

Similarly $B(\widetilde{C}_n)$  corresponds to braids which
are symmetric with respect to rotations about  a vertical line
between $i$ and $i+1$ for $i \equiv 0 \mod n$.
Note that the two rotation symmetries generate a translation,
and thus $B(\widetilde{C}_n)$ is naturally a subgroup of
$B(\widetilde{S}_{2n})$

\begin{figure}[htb]
\begin{center}
\epsfig{figure=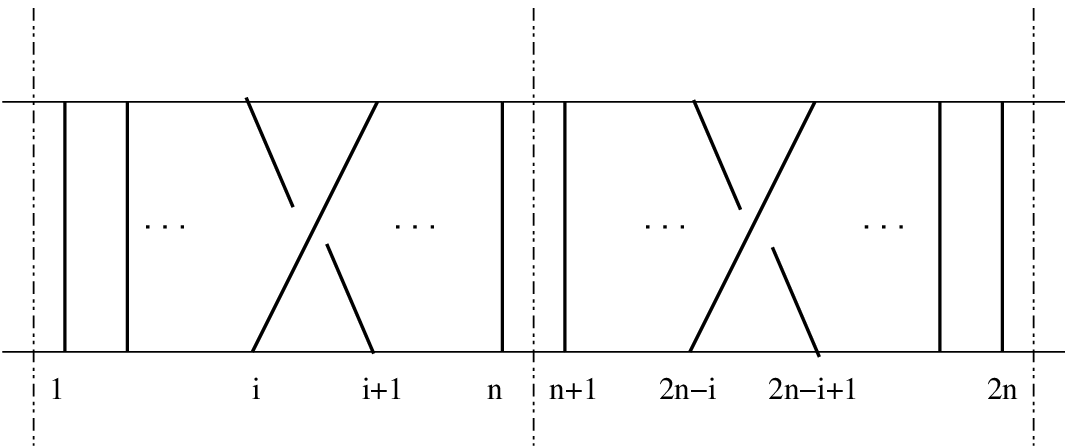}
\end{center}
\caption{$U_i \in B(\widetilde{C}_n)$ }
\end{figure}
\begin{figure}[htb]
\begin{center}
\epsfig{figure=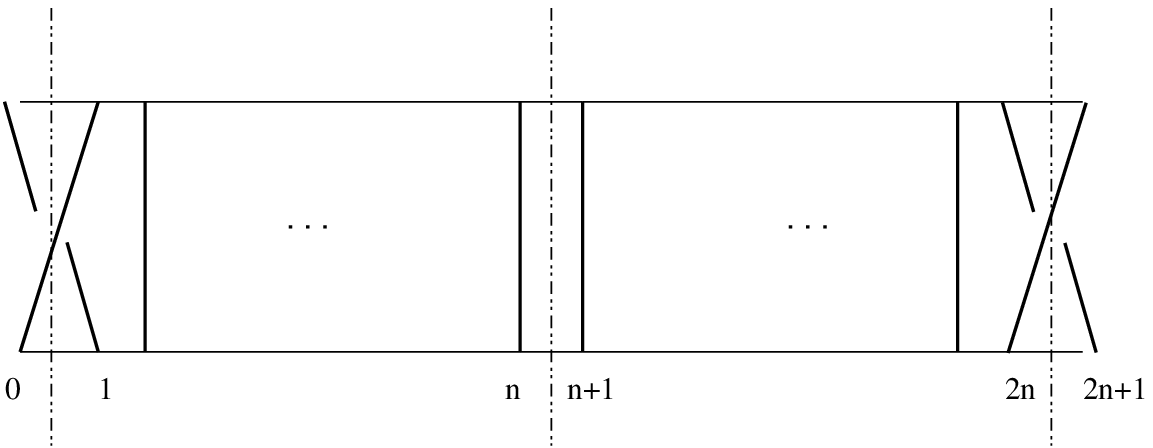}
\end{center}
\caption{$U_0 \in B(\widetilde{C}_n)$ }
\end{figure}
\begin{figure}[htb]
\begin{center}
\epsfig{figure=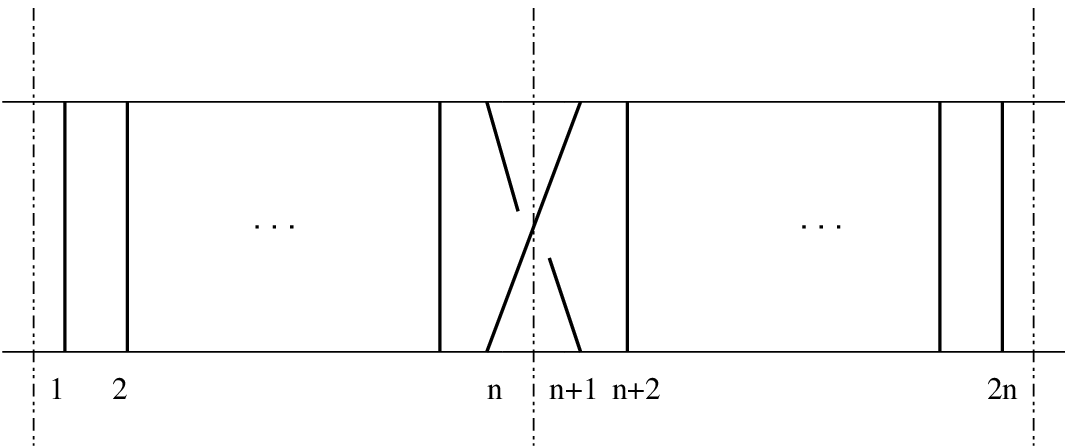}
\end{center}
\caption{$U_n \in B(\widetilde{C}_n)$ }
\end{figure}

\begin{figure}[htb]
\begin{center}
\epsfig{figure=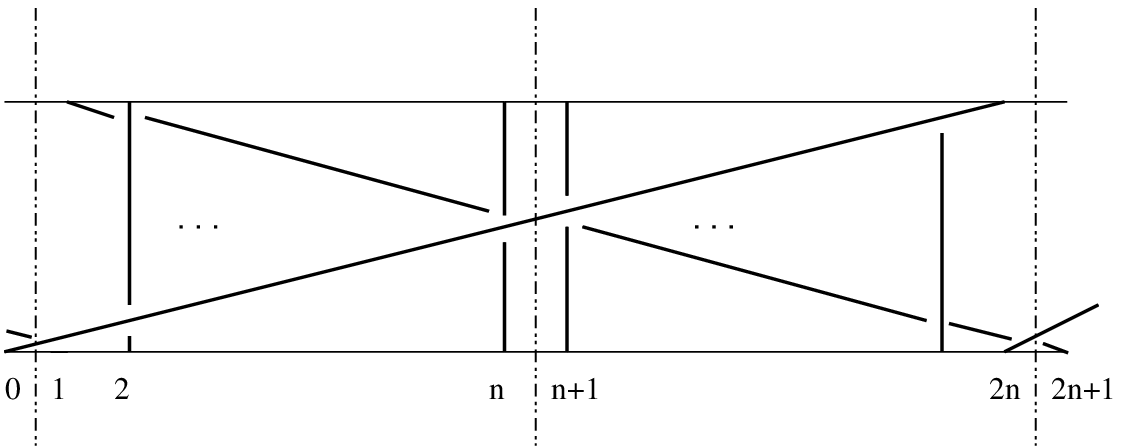}
\end{center}
\caption{$Y_1 \in B(\widetilde{C}_n)$ }
\end{figure}

\begin{figure}[htb]
\begin{center}
\epsfig{figure=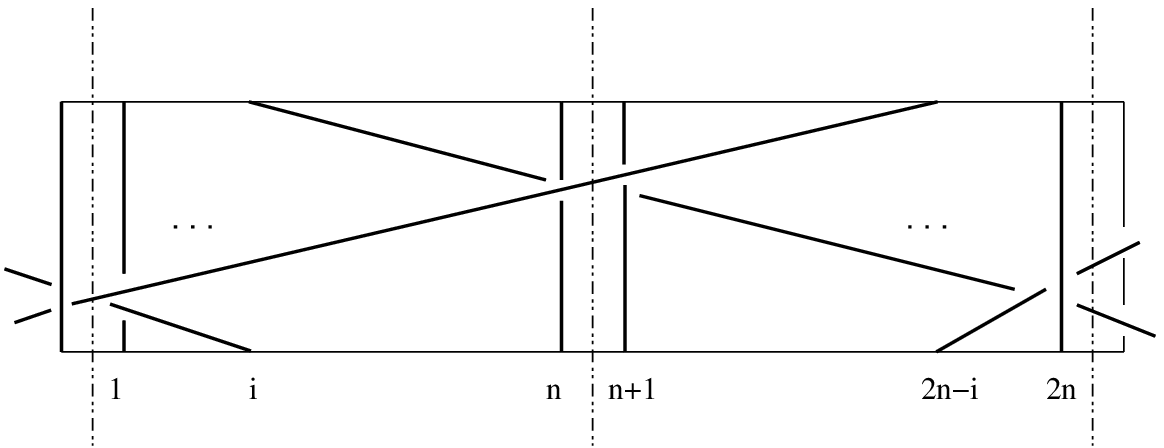}
\end{center}
\caption{$Y_i \in B(\widetilde{C}_n)$ }
\end{figure}

\section{Commutative subgroups of affine braid groups}
\label{sec-comm}

The cleanest proof of our quadratic transformations requires the
construction of nonstandard commutative subalgebras of affine Hecke
algebras.  It turns out that there is a natural construction that
associates a commutative subgroup of an extended 
affine braid group to each chamber
of the associated finite Weyl group.

More precisely, to each chamber we may associate an injective
homomorphism $\Lambda \to B(\widetilde W)$.
We first consider a related construction which associates a map
$\widetilde W\to B(\widetilde W)$ to each alcove of $\widetilde W$.
  For the
standard alcove
this is
just the map 
$$w \mapsto U(w)$$
used to define $B(\widetilde W)$.
More generally we define
\[
\U_{w_1}(w_2):= \U(w_1)^{-1}\U(w_1 w_2).
\]
Note that if we multiply $w_1$ on the right by an element of 
length $0$ that this has no effect on $\U(w_1)$,
which is therefore a function only depending on the associated alcove.
More precisely we have the following.
\begin{lem}
Let $A_0$ denote the standard alcove of the
extended affine Weyl group $\widetilde W$.
  Then for
any simple reflection $s$ of
$\widetilde W$
and any element $w\in \widetilde W$,
\[
\U_w(s) = \U(s)^{\pm 1}
\]
The sign is positive if and only if the simple root corresponding to $s$
is positive for the alcove $A_0^w$.
For any length $0$ element $\sigma$
$$ U_w(\sigma) = \sigma.$$
\end{lem}

\begin{proof}
By definition, we have
\[
\U_w(s) = \U(w)^{-1}\U(w s).
\]
If $\ell(ws)>\ell(w)$, then $\U(ws)=\U(w)\U(s)$, and thus $\U_w(s)=\U(s)$;
otherwise, $\U(w)=\U(ws)\U(s)$, and $\U_w(s)=\U(s)^{-1}$.  Since
$\ell(ws)>\ell(w)$ if and only if $s$ is positive for the alcove $A_0^w$,
the claim follows.

For length 0 elements, we find that $\ell(w\sigma)=\ell(\sigma w)=\ell(w)$,
and the claim follows.
\end{proof}

With this in mind, we will also write
\[
\U_A(w) = \U_{w_A}(w),
\]
where $A$ is the alcove $A_0^{w_A}$.

We also trivially have:

\begin{lem}
For any element $w\in \widetilde W $ and any alcove $A$, 
\[
\U_A(w)^{-1}=\U_{A^w}(w^{-1}).
\]
Similarly, for any elements $w_1$, $w_2\in \widetilde W $
and any alcove $A$, 
\[
\U_A(w_1 w_2) = \U_{A}(w_1)\U_{A^{w_1}}(w_2).
\]
\end{lem}

We can thus describe $\U_A(w)$ as follows: Take any expression
(reduced or not) for $w$ in terms of simple reflections, say
\[
w = s_1 s_2\dots s_n \sigma
\]
Then by iterating the second lemma, we obtain
\begin{align}
\U_A(w)
&=
\U_A(s_1)
\U_{A^{s_1}}(s_2)
\U_{A^{s_1 s_2}}(s_3)
\cdots
\U_{A^{s_1\cdots s_{n-1}}}(s_n) \sigma\\
&=
\U(s_1)^{\pm 1} \U(s_2)^{\pm 1}\cdots \U(s_n)^{\pm 1} \sigma,
\end{align}
where each sign is given by the sign of the given simple root on the
current choice of alcove. 

So far everything we have been saying could apply just as well to
any (extended) Coxeter group.
In the case of an extended affine Weyl group, we have the additional
structure of the associated finite Weyl group $W$.
In particular, in addition to the alcoves of $\widetilde W$,
we may consider the chambers of $W$.

Using the natural quotient map $\widetilde W \to G$ we may associate
to each simple root of $\widetilde W$ a root of $G$ and may
thus sensibly talk about the sign of a root with respect to
 a chamber.  Thus
given a chamber $C$ of the finite Weyl group $W$ and a simple reflection
of $\widetilde W$, we can define
\[
\U_C(s)=\U(s)^{\pm 1},
\]
with positive sign precisely when the corresponding root is positive for
$C$; that is, when the corresponding halfspace contains $C$.  Then for any
word $w=s_ns_{n-1}\dots s_1 \sigma$ in the generators of $\widetilde W$,
we define
\[
\U_C(w)
=
\U_{C}(s_1)\U_{C^{s_1}}(s_2)\U_{C^{s_1s_2}}(s_3)\cdots
\U_{C^{s_1\cdots s_{n-1}}}(s_n)
\sigma.
\]

\begin{thm}
Let $\widetilde W$ be an extended affine Weyl group,
and let $C$ be a chamber of the
associated finite Weyl group $W$.  Then for any word $w$ in the generators of
$\widetilde W$, there exists a vector $v_w$ such that for any alcove
$A\subset v_w+C$,
\[
\U_A(w)=\U_C(w).
\]
In particular, $\U_C(w)$ depends on $w$ only via its image in $\widetilde W$,
and
\[
\U_C(w_1 w_2) = \U_C(w_1)\U_{C^{w_1}}(w_2).
\]
\end{thm}

\begin{proof}
We restrict our attention to the case $\widetilde W = W \ltimes \Lambda_0$;
the general case is analogous.
Write $w=s_1\dots s_n$, and consider
\[
\U_C(w)=\U(s_1)^{\pm 1}\cdots\U(s_n)^{\pm 1}.
\]
For $1\le i\le n$, let $H_i(w)$ denote either the half-space corresponding
to $s_i$ or its complement (the former precisely when $\U(s_i)^{\pm 1}$
occurs with positive sign), and define a sequence of convex sets $D_i(w)$
by:
\begin{align}
D_n(w)&=H_n(w)\\
D_i(w)&=D_{i+1}(w)^{s_i}\cap H_i(w).
\end{align}
We claim that the following is true for
$1\le i\le n$:
\begin{itemize}
\item[(a)] The set $D_i(w)$ is nonempty, and satisfies
\[
D_i(w)+C^{s_1\cdots s_{i-1}}=D_i(w).
\]
\item[(b)] For any alcove $A\subset D_i(w)$,
\[
\U_A(s_i\cdots s_n)
=
\U_{C^{s_1\cdots s_{i-1}}}(s_i\cdots s_n)
\]
\end{itemize}
Indeed, a simple induction argument reduces to the case $n=1$, in which
case $(a)$ and $(b)$ are immediate.

Thus any choice $v_w\in D_1(w)$ proves the first claim of the theorem.
The remaining claims follow from the corresponding results for alcoves.
\end{proof}

The point of using chambers rather than alcoves is that chambers are left
invariant by translations.  As a result, if $\Lambda$ denotes the
translation subgroup of $\widetilde W$,
we find the following.

\begin{cor}
\label{cor-conjugate.braid}
For any chamber $C$, $\U_C$ induces a
homomorphism $\U_C:\Lambda\to B(\widetilde W)$.
The homomorphisms associated to different choices of $C$ are conjugate, in
the sense that
\[
\U_{C^w}(\tau_\nu)
=
\U_C(w)^{-1} \U_C(w \tau_\nu w^{-1}) \U_C(w)
\]
for arbitrary $w\in \widetilde W$.
\end{cor}

\begin{proof}
The first claim is immediate.  For the second claim, we write
\[
\U_C(w\tau_\nu w^{-1})
=
\U_{C}(w)
\U_{C^w}(\tau_\nu) \U_{C^w}(w^{-1}).
\]
\end{proof}
\begin{rem}
Note more generally that for each chamber $C$ we can extend
this homomorphism to a homomorphism from the stabilizer of $C$
to $B(\widetilde W)$.
\end{rem}
We will define $Y^C_\nu =\U_C(\tau_\nu)$ accordingly, and write
$Y_\nu=Y^{C_0}_\nu$.
Note the  $Y^C_\nu$ commute (as $\Lambda$ is commutative).

In addition to the relevance of alternate chambers to our vanishing
results, note
 also that with respect to our standard inner product
for $\widetilde S_n$
 it lets us express the adjoint to the standard $Y_\nu$ as
$Y_{w_0 \nu}^C$ where $C$ is the opposite chamber to the standard one.
\begin{thm}
Suppose the weight $\lambda$ is dominant for the chamber $C$, that is
$\lambda\in C$.  Then
\[
Y^C_\lambda=\U(\tau_\lambda).
\]
In general, if we write $\lambda=\lambda^+-\lambda^-$ with $\lambda^\pm\in
C$, then
\[
Y^C_\lambda=\U(\tau_{\lambda^+})\U(\tau_{\lambda^-})^{-1}.
\]
\end{thm}

\begin{proof}
Let $w$ be a word expressing $\tau_\lambda$ in terms of the generators of
$\widetilde W$,
and choose $v_w$ accordingly.  In particular, we can choose $v_w$ to
be a dominant weight $\lambda'$ for $C$.  We thus find
\[
Y^C_\lambda = \U_{A_0+\lambda'}(\tau_\lambda)
            = \U(\tau_\lambda \tau_{\lambda'})\U(\tau_{\lambda'})^{-1}
\]
But since both $\lambda$ and $\lambda'$ are dominant for $C$, it follows
that
\[
\ell(\tau_\lambda)+\ell(\tau_{\lambda'})=\ell(\tau_{\lambda+\lambda'}),
\]
and thus
\[
\U(\tau_\lambda \tau_{\lambda'})=\U(\tau_\lambda)\U(\tau_{\lambda'});
\]
the result follows.
\end{proof}

In particular, we find that $Y_\lambda$ agrees with the standard
construction of a commutative subgroup of $B(\widetilde W)$.

\begin{thm}
Let $H(\widetilde W)$ be the Hecke algebra corresponding to $\widetilde W$, and let
$Y^C_\nu$ denote the image in $H(\widetilde W)$ of the corresponding element of
$B(\widetilde W)$.  Then for any weight $\lambda \in \Lambda_0$, the sum
\[
\sum_{\mu\in \lambda^{W_0}} Y^C_\mu
\]
is a central element of $H(\widetilde W)$ independent of the choice of chamber $C$.
\end{thm}

\begin{proof}
If we write $C=C_0^w$, then
\[
\sum_{\mu\in \lambda^{W_0}} Y^C_\mu
=
\U(w)^{-1} \sum_{\mu\in \lambda^{W_0}} Y_\mu \U(w),
\]
and thus the claim follows from the standard fact that
\[
\sum_{\mu\in \lambda^{W_0}} Y_\mu
\]
is central.
\end{proof}
\begin{rem}
For general $\lambda \in \Lambda$, this element commutes with all
of the generators, but might act nontrivially on scalars.
\end{rem}
For each $\lambda$ in the root lattice of $\widetilde W$, 
 we can thus define
nonsymmetric Macdonald polynomials $E^C_\lambda$ by
\[
E^C_\lambda \propto \U(w)^{-1} E_{w^{-1}\lambda},
\]
where $C=C_0^w$, and the constant is chosen to make the coefficient of
$x^\lambda$ in $E^C_\lambda$ equal to 1.

%

\bibliographystyle{plain}
\bibliography{quadtrans}

\end{document}